\documentclass[11pt]{amsart}
\usepackage{amsmath,amssymb,latexsym,dsfont}
\usepackage[left=2.6cm,right=2.6cm,top=2.7cm,bottom=2.7cm]{geometry}

\usepackage{amsmath,amsfonts,amssymb,epsfig, textcomp}
\usepackage{graphicx,tikz}
\usepackage{hyperref, cleveref} 
\usepackage{cases,listings}
\usepackage{color}
\usepackage{mathabx}

\newtheorem{theorem}{Theorem}[section]
\newtheorem{lemma}{Lemma}[section]

\newtheorem{proposition}{Proposition}[section]

\newtheorem{remark}{Remark}[section]

\newtheorem{definition}{Definition}[section]
\numberwithin{equation}{section}

      \newcommand{\hy}{\hat y}
      \newcommand{\hu}{\hat u}
      \newcommand{\hv}{\hat v}

      \newcommand{\cH}{\mathcal H}

 \newcommand{\cB}{{\mathcal B}}

   \newcommand{\ty}{\widetilde y}

      \newcommand{\N}{\mathbb{N}}

      \newcommand{\loc}{\operatorname{loc}}
      
      \newcommand{\eps}{\varepsilon}
      \newcommand{\mR}{\mathbb{R}}

      \newcommand{\mC}{\mathbb{C}}

      \newcommand{\dsp}{\displaystyle}
      
\newcommand{\bH}{\mathbf H}

      \makeatletter
      \def\@setcopyright{}
      \def\serieslogo@{}
      \makeatother
   \newcommand{\tr}{^\mathsf{T}}

\newcommand{\hmb}[1]{\textcolor{black}{#1}}
\newcommand{\hmr}[1]{\textcolor{black}{#1}}

\newcommand{\cL}{\mathcal L}
\newcommand{\proj}{\mbox{proj}}

\newcommand{\tu}{\widetilde u}

\newcommand{\cT}{{\mathcal T}}

\newcommand{\cA}{{\mathcal A}}

\newcommand{\cR}{{\mathcal R}}

\newcommand{\tPsi}{{\widetilde \Psi}}

\newcommand{\tz}{\widetilde z}

\newcommand{\be}{\begin{equation}}
\newcommand{\ee}{\end{equation}}

\newcommand{\cF}{{\mathcal F}}

\newcommand{\tf}{{\widetilde f}}

\newcommand{\cD}{{\mathcal D}}

\newcommand{\mH}{\mathbb{H}}

\newcommand{\bR}{{\bf R}}

\newcommand{\wPhi}{\widetilde{\Phi}}
 
\newcommand{\cX}{{\mathcal X}}

\newcommand{\cK}{{\mathcal K}}

\newcommand{\cU}{{\mathcal U}}

\newcommand{\cQ}{{\mathcal Q}}

\newcommand{\bA}{{\bf A}}

\newcommand{\ff}{{\bf f}_{\bH}}

\newcommand{\bh}{{\bf h}}

\newcommand{\bg}{{\bf g}}

\newcommand{\bT}{{\bf T}}

\newcommand{\TT}{{\cT}}

\newcommand{\hlambda}{{\hat \lambda}}

\newcommand{\tv}{\widetilde v}

\newcommand{\bPhi}{{\bf \Phi}}

\newcommand{\ta}{\widetilde \alpha}

\title[Stabilization of bilinear Schr\"odinger equation]{Rapid stabilization and finite time stabilization of the bilinear Schr\"odinger equation}

\author[H.-M. Nguyen]{Hoai-Minh Nguyen}
\address[H.-M. Nguyen]{\hmr{Sorbonne Universit\'e, Universit\'e Paris Cit\'e, CNRS, INRIA, \newline \indent 
Laboratoire Jacques-Louis Lions, LJLL, F-75005 Paris, France}
}
\email{hoai-minh.nguyen@sorbonne-universite.fr}

\begin{document}

\maketitle

\begin{abstract} We propose a method to  establish the rapid stabilization of the bilinear Schr\"odinger control system and its linearized system, and the finite-time stabilization of the linearized system using the Grammian operators. The analysis of the rapid stabilization involves a new quantity (variable) which is inspired by the adjoint state in the optimal control theory and is proposed in our recent work on control systems associated with strongly continuous group, \hmb{and a new regularizing effect of the control system}. The analysis of the finite-time stabilization follows the strategy introduced by Coron and Nguyen in the study of the finite-time stabilization of the heat equation and incorporates a new ingredient involving the estimate of the cost of controls of the linearized system in small time derived in this paper. 
\end{abstract}

\tableofcontents

\section{Introduction}

\subsection{Statement of the main results}
We consider the following bilinear control Schr\"odinger system, with $I = (0, 1)$,  
\be \label{sys-NL}
\left\{\begin{array}{cl}
i \Psi_t = - \Delta \Psi - u(t) \mu(x) \Psi (t, x) & \mbox{ in } \mR_+ \times I, \\[6pt] 
\Psi(t, 0) = \Psi(t, 1) = 0 &  \mbox{ in } \mR_+, \\[6pt]
\Psi(0) = \Psi_0 &  \mbox{ in } I, 
\end{array} \right. 
\ee
where $\Psi_0$ is the initial data, 
$$
\mbox{the control $u$ is {\it real}}, 
$$
and $\mu$ is a given {\it real} function, around the \hmr{ground} state. Here $\Psi$ is the complex-valued wave function of a particle confined in a $1d$ infinite square potential well. The particle is subjected to an electric field inside the domain with the amplitude $u$, and $\mu$ is the dipolar moment of the particle.  
For detailed approximations leading to this first-order interaction Hamiltonian we refer for example to \cite[Chapter 2]{Alessandro}.

Let $\lambda_1 < \lambda_2 < \dots, \lambda_k < \dots $ be the set of eigenvalues of the Laplace equation in $I$ with the zero Dirichlet boundary condition, and let $(\varphi_k)$ be the standard orthogonal basis in $L^2(I)$ formed by the corresponding eigenfunctions. Thus 
$$
\left\{\begin{array}{cl}
-\Delta \varphi_k = \lambda_k \varphi_k & \mbox{ in } I, \\[6pt]
\varphi_k = 0 &  \mbox{ on } \partial I. 
\end{array} \right. 
$$ 
Explicitly, for $k \ge 1$, 
\be \label{def-lambdak-varphik}
\lambda_k = \pi^2 k^2   \quad \mbox{ and } \quad \varphi_k (x) = \sqrt{2} \sin (\pi k x) \mbox{ in } I. 
\ee
It is clear that 
\be \label{property-varphi1}
e^{-i \lambda_1 t} \varphi_1 \mbox{ is a solution of \eqref{sys-NL} with $u = 0$ and $\Psi_0 = \varphi_1$.}
\ee

We are interested in the stabilization of the system \eqref{sys-NL}  around the \hmr{ground} state $e^{- i \lambda_1 t} \varphi_1$. To this end, it is convenient to introduce
\be \label{def-tPsi}
\tPsi (t, x) = e^{i \lambda_1 t} \Psi(t, x) \mbox{ in } \mR_+ \times I \quad \mbox{ and } \quad  \tPsi_0 (x) = \Psi_0 (x) \mbox{ in } I. 
\ee
We then have, by \eqref{sys-NL},  
\be \label{sys-NL2}
\left\{\begin{array}{cl}
i \tPsi_t = - \Delta \tPsi  - \lambda_1 \tPsi - u(t) \mu(x) \tPsi (t, x) & \mbox{ in } \mR_+ \times I, \\[6pt] 
\tPsi(t, 0) = \tPsi(t, 1) = 0 &  \mbox{ in } \mR_+. 
\end{array} \right. 
\ee
The linearized system of \eqref{sys-NL2} when $\tPsi$ is closed to $\varphi_1$, i.e., $\Psi(t, x)$ is closed to $e^{-i \lambda_1 t} \varphi_1$, is 
\be\label{sys-LN2}
\left\{\begin{array}{cl}
i \tPsi_t = - \Delta \tPsi  - \lambda_1 \tPsi - u(t) \mu(x) \varphi_1 (x) & \mbox{ in } \mR_+ \times I, \\[6pt] 
\tPsi(t, 0) = \tPsi(t, 1) = 0 &  \mbox{ in } \mR_+. 
\end{array} \right. 
\ee 

In what follows,  we {\it always} assume that  
\be
\mu \in H^3(I, \mR).
\ee
The following condition on $\mu$ is used later:  
\be \label{cond-mu}
|\langle \mu \varphi_1,  \varphi_k \rangle_{L^2(I)}| \ge \frac{c}{k^3} \mbox{ for } k \in \N_+,  
\ee
for some positive constant $c$ unless stated differently  \footnote{Hereafter, given a Hilbert space $\cH$, we denote $\langle \cdot, \cdot \rangle_{\cH}$ its scalar product.} 
We are interested in the solutions of the above Schr\"odinger systems with controls $u$ in $L^2_{loc}([0, + \infty); \mR)$ (we insist again that we are interested in the controls which are {\it real}). 

The condition \eqref{cond-mu} is a sufficient condition to have the exact controllability of the linearized systems in small time and this implies the local exact controllability of the nonlinear systems, as shown by Beauchard and Laurent \cite{BL10}. This condition is also a necessary condition to ensure that the nonlinear systems are locally exactly controllable in small time, see the work of Beauchard and Morancey \cite{BM14}. The condition \eqref{cond-mu} is generic, see \cite[Appendix A]{BL10}. 

\medskip 

\hmr{In this paper, we first present a method to obtain the rapid stabilization of the linearized control system \eqref{sys-LN2}, and of the bilinear control system \eqref{sys-NL2} locally around the ground state. This involves the Gramian operators, see \eqref{def-Q}, and is done under the exact controllability assumption of the linearized system \eqref{sys-LN2} given in \eqref{cond-mu-EEE}; this is weaker than \eqref{cond-mu}. More precisely, concerning \eqref{sys-LN2}, under the exact controllability assumption, we will show that, for every $\lambda > 0$, there exists a linear operator (feedback) $K$ (involving the Gramian operators and depending on $\lambda$, see \eqref{def-uuu}), appropriately defined and taking {\it real} values, such that, for some positive constant $C$, it holds  
$$
\| \tPsi (t, \cdot)\| \le C e^{- \lambda t} \| \tPsi (0, \cdot)\|, 
$$
for all initial data $\tPsi (0, \cdot)$ in an appropriate space,  where $\tPsi$ is the solution of \eqref{sys-LN2} with $u(t) = K \tPsi(t, \cdot)$. Similar facts hold for \eqref{sys-NL2} for all initial data sufficiently closed to $\varphi_1$ in a suitable norm. We then extend this approach to obtain the finite-time stabilization of the linearized control system \eqref{sys-LN2} using a piecewise constant feedback under the assumption \eqref{cond-mu}.}


\medskip
As in previous works, see, e.g., \cite{BL10, BM14,CGM18}, we are interested in the solutions in the space $\bH$ (for each time $t$) defined by 
\be \label{def-bH}
\bH = \Big\{ \Psi \in H^1_0(I; \mC);  \sum_{k \ge 1} |k^3 \langle \Psi, \varphi_k \rangle_{L^2(I)}|^2 < + \infty \Big\},  
\ee
and 
\be \label{def-bH1}
\bH_{1, \sharp} = \Big\{ \Psi \in \bH \mbox{ such that }  \Re \langle \Psi, \varphi_1 \rangle_{L^2(I)} = 0 \Big\}   
\ee
(\hmr{a property related to \eqref{def-bH1} is given in \eqref{contraint}}). Here and in what follows, for a complex number $z$, we denote its real part, its complex part, and its complex conjugate by $\Re z$, $\Im z$, and $\bar z$, respectively. 
We equip the following scalar product for the spaces  $\bH$ and $\bH_{1, \sharp}$: 
\be \label{scalar-product-H}
\langle \Psi, 
\tPsi \rangle_{\bH}  = \langle \Psi, \tPsi \rangle_{H^3(I)} := \int_I 
 \Big( \Psi \overline{\tPsi} + \Psi'  \overline{\tPsi'} + \Psi'' \overline{\tPsi}'' + \Psi''' \overline{\tPsi}''' \Big) \, ds \mbox{ for } y, \ty \in \bH,  
\ee
and 
\be \label{scalar-product-H1}
\langle \Psi, 
\tPsi \rangle_{\bH_{1, \sharp}} = \langle y, 
\ty \rangle_{\bH}   \mbox{ for } \Psi, \tPsi \in \bH_{1, \sharp}.   
\ee

One can show that, for the linearized system \eqref{sys-LN2}, 
\be\label{contraint}
\Psi (t) \in \bH_{1, \sharp} \mbox{ for } t \ge 0 \mbox{ if } \Psi_0 \in \bH_{1, \sharp}. 
\ee
This property does {\it not} hold for the nonlinear system.  Note that the exact controllability has been established for solutions in $C([0, T]; \bH)$, which requires roughly three derivatives in the space variable of the solutions. It is known from a general result of Ball, Marsden, and Slemrod \cite{BMS82} that the Schr\"odinger system \eqref{sys-NL2} is not exactly controllable for solutions in $C([0, T]; H^1_0(I))$ or in $C([0, T]; H^1_0(I) \cap H^2(I))$ when $\mu$ is smooth since the control operator is bounded in this case. 

\hmr{Since we deal only with real controls}, it is convenient to consider the real part and the imaginary part of $\tPsi$ separately. Assume that 
$$
\tPsi = \tPsi_1 + i \tPsi_2 \mbox{ in } \mR_+ \times I, 
$$
where $\tPsi_1$ and $\tPsi_2$ are the real and the imaginary parts of $\tPsi$. System \eqref{sys-NL2} can be written under the form 
\be\label{sys-NL3}
\left\{\begin{array}{cl}
\tPsi_{1, t} = -  \Delta \tPsi_2 - \lambda_1 \tPsi_2 - u(t) \mu(x) \tPsi_{2} & \mbox{ in } \mR_+ \times I, \\[6pt]
\tPsi_{2, t} = \Delta \tPsi_1 + \lambda_1 \tPsi_1 + u(t) \mu(x) \tPsi_1 & \mbox{ in } \mR_+ \times I,  
\end{array} \right.
\ee
and system \eqref{sys-LN2} can be written under the form 
\be\label{sys-LN3}
\left\{\begin{array}{cl}
\tPsi_{1, t} = -  \Delta \tPsi_2 - \lambda_1 \tPsi_2 & \mbox{ in } \mR_+ \times I, \\[6pt]
\tPsi_{2, t} = \Delta \tPsi_1 + \lambda_1 \tPsi_1 + u(t) \mu(x) \varphi_1 & \mbox{ in } \mR_+ \times I. 
\end{array} \right.
\ee

Denote
\be
\mH = \Big\{ y= (y_1, y_2)\tr \in H^1_0(I; \mR^2);  \sum_{\ell=1}^2 \sum_{k \ge 1} |k^3 \langle y_\ell, \varphi_k \rangle_{L^2(I)}|^2 < + \infty \Big\},  
\ee
and 
\be
\mH_{1, \sharp} = \Big\{ y = (y_1, y_2)\tr \in \mH \mbox{ such that }  \langle y_1, \varphi_1 \rangle_{L^2(I)} = 0 \Big\},  
\ee
and we equip with the following scalar products for the spaces  $\mH$ and $\mH_{1, \sharp}$: 
\begin{multline} \label{scalar-product-H}
\langle y, 
\ty \rangle_{\mH}: =  \langle y, 
\ty \rangle_{H^3(I)} = \int_I 
\sum_{\ell =1}^2  \Big( y_{\ell} \ty_{\ell} + y_{\ell}'  \ty_{\ell}' + y_{\ell}'' \ty_{\ell}'' + y_{\ell}''' \ty_{\ell}''' \Big) \, ds \\[6pt] 
\mbox{ for } y = (y_1, y_2)
\tr, \ty = (\ty_1, \ty_2)\tr \in \mH,  
\end{multline}
and 
\be \label{scalar-product-H1}
\langle y, 
\ty \rangle_{\mH_{1, \sharp}} = \langle y, 
\ty \rangle_{\mH}   \mbox{ for } y, \ty \in \mH_{1, \sharp}.   
\ee

It is clear that 
$$
\Psi \in \bH \mbox{ if and only if } (\Psi_1, \Psi_2)\tr \in \mH
$$
and 
$$
\Psi \in \bH_{1, \sharp} \mbox{ if and only if } (\Psi_1, \Psi_2)\tr \in \mH_{1, \sharp},  
$$
where $\Psi_1$ and $\Psi_2$ are the real part and the imaginary part of $\Psi$, respectively. 

One can check, see e.g., \cite{CGM18}, that 
\be \label{CGM1}
\mH = \Big\{ y = (y_1, y_2)\tr \in H^3(I; \mR^2);  y_1(x) =y_2(x) = y_1''(x) =  y_2''(x) = 0 \mbox{ on } \partial I \Big\}
\ee
and 
\begin{multline} \label{CGM2}
\mH_{1, \sharp} = \Big\{ y = (y_1, y_2)\tr \in H^3(I; \mR^2);  \\[6pt] y_1(x) =y_2(x) = y_1''(x) =  y_2''(x) = 0 \mbox{ on } \partial I  \mbox{ and }  \langle y_1, \varphi_1 \rangle_{L^2(I)} = 0 \Big\}.
\end{multline}
Note that $\bH_{1, \sharp}$ is not a subspace of $\bH$ (with respect to the scalar field $\mC$) whilst $\mH_{1, \sharp}$ is a subspace of $\mH$ (with respect to the scalar field $\mR$).

Consider $A: \cD(A) \subset \mH \to \mH$ defined by 
\be \label{def-A}
A y = 
\left(\begin{array}{cc} - \Delta y_2 - \lambda_1 y_2 \\[6pt]
\Delta y_1 + \lambda_1 y_1
\end{array}\right)
 \quad \mbox{ and } \quad \cD(A) = \Big\{ y \in \mH; A y \in \mH \Big\}.  
\ee
Then $\cD(A)$ is dense in $\mH$ and $A$ is skew-adjoint (see \Cref{lem-density,lem-A}). We equip $\cD(A)$ with the standard scalar product for the graph-norm and denote $\cD(A)'$ the dual space of $\cD(A)$. It is worth noting that \hmr{the} definition of $A$ and the domain $\cD(A)$ \hmr{given here} are different from \cite{BL10, CGM18}. Our definitions are motivated by the theory of stabilization developed for control systems associated with a strongly continuous group \cite{Ng-Riccati} and will be clear later when the feedback operator is introduced (see, e.g., \eqref{def-Q}, see also \eqref{def-B*}). 

Let $(A^*, \cD(A^*))$ denote the adjoint of $(A, \cD(A))$, and let 
$$
\mbox{$B : \mR \to \cD(A^*)'$}
$$  be defined by, with $y = (y_1, y_2)\tr \in \cD(A^*)$,  
\be \label{def-B}
\langle Bu, y \rangle_{\cD(A^*)', \cD(A^*)} 
\\[6pt]
= u \Big( \langle \mu \varphi_1, y_2 \rangle_{H^3(I)}  -  (\mu \varphi_1)_{xx}(1)  y_{2, xxx}(1)
+  (\mu \varphi_1)_{xx}(0)  y_{2, xxx}(0)
 \Big). 
\ee
The linear system \eqref{sys-LN3} can be written under the form 
\be \label{sys-LN1-AB}
y' = A y + B u \mbox{ in } \mR_+,
\ee
and the nonlinear system \eqref{sys-NL3} can be written under the form 
\be \label{sys-NL1-AB}
y' = A y + B u + u F (y - \Phi_1) \mbox{ in } \mR_+,
\ee
where 
\be 
\Phi_1 = (\varphi_1, 0)\tr \quad \mbox{ and } \quad \label{def-F}
F(y) = (-  \mu y_2,  \mu y_1)\tr, 
\ee
and, for all $\varphi \in \cD(A^*)$,  
\begin{multline} \label{def-uF}
\langle u F(y), \varphi \rangle_{\cD(A^*)', \cD(A^*)} \\[6pt]
= u \Big( \langle F(y), \varphi \rangle_{H^3(I)}  - \langle (F(y))_{xx}(1),  \varphi_{xxx}(1) \rangle_{\mR^2} 
+ \langle (F(y))_{xx} (0),  \varphi_{xxx} (0) \rangle_{\mR^2} 
\Big)
\end{multline}
(see \Cref{lem-product,lem-WP-A}). 

One cannot extend  $B$ as a bounded operator from $\mR$ into $\mH$. This is the main source of the difficulties in the study of the stabilization using feedback of the linearized system \eqref{sys-LN1-AB} and more critical in the study of the nonlinear system \eqref{sys-NL1-AB} since $u F (y) \not \in L^1((0, T); \mH)$. Nevertheless,  $B$ is an {\it admissible} control operator with respect to the semi-group $\big( e^{t A}\big)_{t \ge 0}$ generated by $A$ in the sense that, for all $u \in L^2([0, T]; \mR)$, it holds that 
\be \label{cond-admissibility}
y \in C([0, T]; \mH) \mbox{ where } y(t): = \int_0^t e^{(t-s) A} B u(s) \, ds 
\ee
(see \Cref{lem-WP,lem-product}).  As a consequence of the closed graph theorem, see e.g., \cite{Brezis-FA}, one has 
\be \label{cond-admissibility-a}
\| y \|_{C([0, T]; \mH)} \le C_T \| u \|_{L^2((0, T); \mR)}. 
\ee
Thus, see e.g., \cite{Coron06,TW09},  that, for $T > 0$, there exists $C_T > 0$ such that 
\hmr{\be \label{ineq-direct-***}
\int_0^T | B^* e^{ sA^*} z|^2 \hmr{\, ds}  \le C_T \|z\|_{\mH}^2 \mbox{ for all } z \in \mH,   
\ee
where $B^*$ is  the adjoint of $B$, and $(e^{t A^*})_{t \in \mR}$ is the group generated by $A^*$. 
This implies 
\be 
\int_0^T | B^* e^{ (T-s) A^*} z|^2 \hmr{\, ds}  \le C_T \|z\|_{\mH}^2 \mbox{ for all } z \in \mH,   
\ee
Since $A$ is skew-adjoint, $e^{ (T-s) A^*} z = e^{-s A^*} e^{T A^*}z$, it follows that}
\be \label{ineq-direct}
\int_0^T | B^* e^{ -sA^*} z|^2 \hmr{\, ds}  \le C_T \|z\|_{\mH}^2 \mbox{ for all } z \in \mH.
\ee

Note that 
$$
B^* :  \cD(A^*) \to \mR. 
$$
and, with $v = (v_1, v_2)\tr \in  \cD(A^*)$, which is also $\cD(A)$ since $A$ is skew-adjoint,  
\be \label{def-B*}
B^* v  =  \langle \mu \varphi_1,  v_2 \rangle_{H^3(I)}  - (\mu \varphi_1)_{xx} (1)  v_{2, xxx}(1) 
+ (\mu \varphi_1)_{xx}(0)  v_{2, xxx}(0) 
\ee
since 
\be
\langle B u, v \rangle_{\mH} = \langle u, B^*v \rangle_{\mR}.
\ee

This paper is devoted to the stabilization of the nonlinear system \eqref{sys-NL2} and its linearized system \eqref{sys-LN2}. 
The rapid stabilization of the linearized control system \eqref{sys-LN2} \hmr{under the assumption \eqref{cond-mu}} was established by Coron, Gagnon, and Morancey \cite{CGM18} using techniques related to backstepping methods. The idea is to transform the original system into a damping one for which the stabilization is an easier task. Their transformations are of Fredholm type and different from the standard Volterra ones in the backstepping method. The existence of these transformations is ensured by \eqref{cond-mu}. The main technical difficulty in the work of Coron, Gagnon, and Morancey \cite{CGM18} is to deal with a control operator that is only admissible but {\it not} bounded. It is worth noting that the backstepping technique and its extended versions are useful tools to stabilize various equations in one-dimensional space such as heat equations \cite{Liu00}, Schr\"odinger equations \cite{KGS11}, KdV equations \cite{CC13,CL14}, hyperbolic systems \cite{CVKB13,CoronNg19,CoronNg19-2} and the references therein. The backstepping can be also used to get finite-time stabilization for heat equations, see \cite{CoronNg17}.  A concise introduction to the backstepping technique can be found in \cite{Krstic08}. At this stage, to our knowledge, \cite{CGM18} is the only work dealing with the rapid stabilization of the linearized Schr\"odinger system using bilinear controls, and the analysis in \cite{CGM18} has not been successfully extended to the nonlinear system.

The goal of this paper is to present another method to obtain the rapid stabilization of the linearized control system \eqref{sys-LN2}  and of the bilinear control system \eqref{sys-NL2} \hmr{under the exact controllability assumption of the linearized system \eqref{sys-LN2}}, and the finite-time stabilization of the linearized control system \eqref{sys-LN2} \hmr{under the assumption \eqref{cond-mu}}. 
Our approach is inspired by our recent work  \cite{Ng-Riccati} in which we study the stabilization of systems associated with a strongly continuous group for unbounded control operators using Gramian operators. \hmb{The new regularizing effect of the control system (see \Cref{lem-WP}) reveals in this paper will play an important role in our analysis. In order to apply this, we also need to derive new information of the Gramian operator, and this is given in \Cref{sect-Q}.} For control systems associated with a strongly continuous group, under the assumption that the systems are exactly controllable, it is shown in \cite{Ng-Riccati} that one can obtain rapid stabilization using static feedback in a trajectory sense or using dynamic feedback. The static trajectory feedback has its roots in the linear quadratic optimal control theory, as developed in Flandoli, Lasiecka, and Triggiani \cite{FLT88} (see also \cite{LT91,WR00,Zwart96,Staffans05,TWX20}). It is known from the optimal control theory that there exists static feedback in a {\it weak} sense to rapidly stabilize the system. Such feedback is understood in a weak sense 
since it is defined only on a dense set of the space state depending on the feedback operator (see \cite[Proposition 4.1]{Ng-Riccati}). The use of Gramian operators to rapidly stabilize exactly controllable systems associated with a strongly continuous group has been previously considered in \cite{Komornik97,Urquiza05,Vest13} via the optimal control theory, and the feedback is thus understood in the {\it weak} sense. One cannot use Gramian operators to stabilize nonlinear settings using the theories developed in \cite{Komornik97,Urquiza05,Vest13} as discussed in \cite{Ng-Riccati} (see also \cite{CGM18}).  In this paper, we show that, for the considered bilinear control Schr\"odinger systems, even if the control operator is unbounded one can still obtain static feedback. Moreover, we construct piecewise constant feedback to reach the stabilization in finite-time for the linearized system.

Before introducing the feedback, we state the observability inequality for the exact controllability of the linearized system \eqref{sys-LN2} in time $T_0$.

\begin{lemma}\label{lem-Obs} Let $\mu \in H^3(I; \mR)$. Assume that
\be \label{cond-mu-EEE}
\mbox{ the linearized system \eqref{sys-LN2} is exactly controllable in $\mH_{1, \sharp}$ in time $T_0$,}
\ee 
for some $T_0 > 0$. We have 
\be \label{cond-mu-EEE-cl}
\int_0^{\hmr{T_0}} | B^* e^{-s A^*} z|^2  \ge C_{T_0} \| z\|_{\mH}^2 \mbox{ for all } z \in \mH_{1, \sharp}, 
\ee
for some positive constant $C_{T_0}$ independent of $z$. 
\end{lemma}

\hmr{Assertion \eqref{cond-mu-EEE-cl} follows from the fact that
\be
\int_0^{\hmr{T_0}} | B^* e^{s A^*} z|^2  \ge C_{T_0} \| z\|_{\mH}^2 \mbox{ for all } z \in \mH_{1, \sharp}, 
\ee
and $A$ is skew-adjoint by the arguments used to obtain \eqref{ineq-direct} from \eqref{ineq-direct-***}.}

\begin{remark} \rm \hmr{As a consequence of the exact controllability result of Beauchard and Laurent \cite{BL10}, and  \Cref{lem-WP,lem-product,lem-WP-A} (see also \Cref{pro-cost} for a more quantitative version for small $T_0$),  the linearized system \eqref{sys-LN2} is exactly controllable in time $T_0$ for all $T_0 > 0$ if \eqref{cond-mu} holds. }
\end{remark}

We are ready to introduce the Gramian operator to stabilize the linearized system \eqref{sys-LN2} and the nonlinear system \eqref{sys-NL2}.  Let $\mu \in H^3(I; \mR)$  and $\lambda > 0$. Define $Q = Q(\lambda): \mH \to \mH$, \hmr{the Gramian operator},  by 
\be \label{def-Q}
\langle Q z, 
\tz \rangle_{\mH} = \int_0^\infty e^{-2 \lambda s} \langle B^* e^{-s A^*} z, B^* e^{-s A^*} \tz \rangle_{\mR} \, ds \mbox{ for } z, \tz \in \mH.
\ee
Since $A$ is skew-adjoint by \Cref{lem-A} in \Cref{sect-pre}, it follows from \eqref{ineq-direct} that $Q$ is well-defined and is symmetric in $\mH$. We also have, by \cite[Proposition 5.1]{Ng-Riccati} (see also \cite{Komornik97,Urquiza05}), 
\be \label{identity-Op-Q}
A Q + Q A^* - B B^* + 2 \lambda Q = 0  
\ee
in the following sense
\be \label{identity-Op-Q-meaning}
\langle Q z, A^* \tz \rangle_{\mH} + \langle A^* z, Q \tz \rangle_{\mH} - \langle B^* z, B^* \tz \rangle_{\mH} + 2 \lambda \langle Q z, \tz \rangle_{\mH} \hmr{=0} \mbox{ for } z, \tz \in \cD(A^*).
\ee
Moreover, if the condition \hmr{\eqref{cond-mu-EEE} holds }then, by \Cref{lem-Obs},
\be \label{def-cQ}
\cQ: = \proj_{\mH_{1, \sharp}} \circ Q : \mH_{1,\sharp} \to \mH_{1, \sharp} \mbox{ is positive, i.e., } 
\langle \cQ z, z \rangle_{\mH_{1, \sharp}} \ge C \| z\|_{\mH_{1, \sharp}} \mbox{ for all } z \in \mH_{1, \sharp}.
\ee
 
\hmb{As shown later in \Cref{lem-Q3}, the application $Q$ maps $\mH_{1, \sharp}$ into $\mH_{1, \sharp}$ and thus the projection $\proj_{\mH_{1, \sharp}}$ can be ignored in the definition of $\cQ$.} 

\medskip

Concerning the rapid stabilization of \eqref{sys-LN3}, we prove the following result. 

\begin{theorem} \label{thm1-S} Let $\mu \in H^3(I, \mR)$ 
and $\lambda  > 0$. \hmr{Assume that   \eqref{cond-mu-EEE} holds}. Given $y_0 \in \mH_{1, \sharp}$, let $y  \in C([0, + \infty); \mH) $ be the unique weak solution of the system 
\be \label{thm1-S-sys}
\left\{\begin{array}{cl}
y' = A y + B u \mbox{ in } \mR_+, \\[6pt]
y(0) = y_0, 
\end{array} \right.
\ee
with \hmr{$u \in L^2_{\loc}([0, +\infty); \mR)$} verifying  
\be \label{def-uuu}
u = -  B^*  \cQ^{-1} \proj_{\mH_{1, \sharp}} y. 
\ee
Then 
\be\label{thm1-S-cl0}
\hmr{y(t) \in \mH_{1, \sharp} \mbox{ for } t \ge 0,}
\ee 
and
\be \label{thm1-S-cl2}
\|\cQ^{-1} y(t) \|_{\mH} = e^{-2 \lambda \hmr{t}} \| \cQ^{-1}y_0 \|_{\mH} \mbox{ for }  t \ge 0. 
\ee
Consequently, there exist two positive constants $C_1, C_2$ independent of $y_0$ such that  
\be \label{thm1-S-decay}
C_1 e^{-2 \lambda t}\| y_0 \|_{\mH}  \le \| y(t)\|_{\mH} \le C_2 e^{-2 \lambda t}\| y_0 \|_{\mH} \mbox{ for }  t \ge 0.  
\ee
\end{theorem}

\begin{remark} \rm The meaning of the weak solutions are given in \Cref{def-WS-cA} in \Cref{sect-appendixA} for which one considers $Bu$ as a source term.  The well-posedness of \eqref{thm1-S-sys} is a part of the conclusion of \Cref{thm1-S}. 
\end{remark}

As a consequence of \Cref{thm1-S}, the linearized system \eqref{sys-LN3} is rapidly stabilizable by feedback controls. Equivalently, the linearized system \eqref{sys-LN2} is rapidly stabilizable by feedback controls. 

\begin{remark} \label{rem-thm1-S} \rm Note that $y(t) \in \mH_{1, \sharp}$ for $t \ge 0$. One can hence replace the term $\proj_{\mH_{1, \sharp}} y$ by $y$ in \Cref{thm1-S}. 
\end{remark}

Concerning the non-linear system \eqref{sys-NL3}, which is equivalent to system \eqref{sys-NL2},  we have the following result. 

\begin{theorem} \label{thm2-S} Let $\mu \in H^3(I, \mR)$. \hmr{Assume that  \eqref{cond-mu-EEE} holds} and let $\lambda > 0$. For $0 < \hlambda < \lambda$, there exist two positive constants $\eps_0 > 0$ and  $C > 0$ such that


\be
\| y(t, \cdot) - \Phi_1 \|_{\mH} \le C  e^{-2\hlambda t} \| y(0, \cdot) -  \Phi_1 \|_{\mH} \quad \mbox{ for } t \ge 0,  
\ee
for all $y_0 \in \mH_{1, \sharp}$ with $\| y_0\|_{L^2(I)} = 1 $ and $\|y_0 -  \Phi_1 \|_{\mH} \le \eps_0$, where $y(t, \cdot) \in C\big([0, T]; \mH \big)$ is the unique weak solution of the system \be \label{thm2-S-sys}
\left\{\begin{array}{cl}
y' = A y + B u + u F(y -  \Phi_1)  \mbox{ in } \mR_+, \\[6pt]
y(0) = y_0, 
\end{array} \right.
\ee
with \hmr{$u \in L^2_{\loc}([0, +\infty); \mR)$} verifying  
$$
u  = - B^* \cQ^{-1} \proj_{\mH_{1, \sharp}} (y -  \Phi_1) = \hmr{- B^* \cQ^{-1} \proj_{\mH_{1, \sharp}} y}.
$$

\end{theorem}

\begin{remark} \rm The meaning of the weak solutions are given in \Cref{def-WS-cA} in \Cref{sect-appendixA} for which one considers $u F(y -  \Phi_1)$ as the source term. 
\end{remark}

As a consequence of \Cref{thm2-S}, the nonlinear bilinear control system \eqref{sys-NL3} is locally rapidly stabilizable by feedback controls. Equivalently, the nonlinear bilinear control system \eqref{sys-LN2} is locally rapidly stabilizable by feedback controls.

\medskip 

Concerning the finite-time stabilization, we have the following result on the linearized bilinear control system \eqref{sys-LN1-AB}.

\begin{theorem} \label{thm-FT-LN}  Let $\mu \in H^3(I, \mR)$ be such that \eqref{cond-mu} holds and let $T > 0$. \hmb{There exists an increasing sequence $(t_n)_{n \ge 0}$ converging to $T$ with $t_0 = 0$, and an increasing sequence of positive numbers $(\lambda_n)_{n \ge 0}$ converging to infinity such that}, for every $y_0 \in \mH_{1, \sharp}$, there exists a unique solution $y \in C([0, T); \mH)$ of the system 
\be
\left\{\begin{array}{cl}
y' = A y + B u \mbox{ in } [0, T), \\[6pt]
y(0) = y_0, 
\end{array} \right.
\ee
with, \hmr{$u \in L^2([0, T); \mR)$} verifying, with   $\cQ_n = \cQ(\lambda_n)$, 
\be
u(t) = -  B^*  \cQ_n^{-1} \proj_{\mH_{1, \sharp}} y (t) \mbox{ in } [t_n, t_{n+1}).  
\ee
Moreover, 
\be
y(t, \cdot) \to 0 \mbox{ in } \mH \mbox{ as } t \to T_{-}
\ee
and 
\be
u(t, \cdot) \to 0 \mbox{ as } t \to T_{-}. 
\ee
\end{theorem}

\begin{remark} \rm \hmr{It is worth noting that for the rapid stabilization, one requires the exact controllability of the linearized system instead of \eqref{cond-mu}. Concerning the finite-time stabilization, \eqref{cond-mu} is used to estimate the control cost for small time. This is employed to derive upper bounds for the norm $\cQ^{-1} = \cQ^{-1} (\lambda)$ where the dependence on $\lambda$ is explicit.}
\end{remark}

\subsection{Ideas of the proof}  The approach used in this paper is inspired by our recent work \cite{Ng-Riccati}. 
We first discuss the analysis of \Cref{thm1-S} and \Cref{thm2-S}. 
Concerning the linearized system \eqref{sys-LN2} (whose results are given in \Cref{thm1-S}), one of the main parts of the analysis is to develop the theory in \cite{Ng-Riccati} to take into account the intrinsic constraint \eqref{contraint}. \hmb{The other important part is to derive new regularizing effect of the linearized control system (see \Cref{lem-WP}) and the new information of the Gramian operator (see \Cref{sect-Q}). These facts and their analysis are interesting in themselves.} 
Concerning the rapid stabilization of the nonlinear system \eqref{sys-NL2}  (whose results are given in \Cref{thm2-S}), in addition to the ingredients used for the linearized system, we essentially use the fact that the solutions of the Schr\"odinger system conserve the $L^2$-norm. This fact is used to control the component of the solution which is orthogonal to $\mH_{1, \sharp}$ with respect to the $L^2(I)$-scalar product (or the $\mH$-scalar product). Additional technical ingredients for all the stabilization results are the well-posedness of the nonlinear feedback control systems, which are of nonlinear, nonlocal, and non-bounded nature (see \Cref{sect-WP}), and the way to translate the results between the original systems \eqref{sys-NL2} and \eqref{sys-LN2} and the corresponding systems written under the semi-group language \eqref{sys-LN1-AB} and \eqref{sys-NL1-AB} (see \Cref{lem-WP-A} and \Cref{sect-appendixA}).  

To take one step further from the rapid stabilization to obtain the finite-time stabilization (\Cref{thm-FT-LN}) for the linearized system, we follow the strategy of Coron and Nguyen \cite{CoronNg17}. The idea is to stabilize the system more and more as the time $t$ goes to $T_-$. To be able to apply the strategy in \cite{CoronNg17}, one needs to understand the size of $\|\cQ_n^{-1}\|_{\cL(\mH_{1, \sharp})}$  \footnote{Hereafter, given a Hilbert space $\cH$, we denote $\cL(\cH)$ the space of all continuous linear applications from $\cH$ to $\cH$ equipped with the standard norm.} as a function of $\lambda_n$ (a good bound for the size of $\|Q_n \|_{\cL(\mH)}$ follows from the admissibility of the control operator $B$, see \eqref{ineq-direct}). This is given in \Cref{lem-Qlambda} after establishing the cost of the control for small time (see \Cref{pro-cost}). This result is interesting in itself and its proof uses similar techniques as in \cite{TT07}. The way to gain suitable information in each time interval $[t_n, t_{n+1})$ here is different from the one in \cite{CoronNg17} for which precise estimates of kernels of transformations from the backstepping technique are derived using the information of the kernels. Our new way to get appropriate information to be able to apply the strategy in \cite{CoronNg17}  is quite robust and can be used in different contexts where the size of the control cost is understood for small time. An application of this approach will be given in \cite{Ng-KdV-stabilization} to study the finite-time stabilization of a KdV control system.

\subsection{Previous related results}

The controllability properties for the Schrödinger equation were mostly studied in the usual linear setting (in contrast to the bilinear control problems considered here). For the control of the linear Schrödinger equation with internal control (localized on a subdomain), we refer to \cite{Lebeau92,Machtyngier}, the survey \cite{Laurent14}, and the references therein. In this setting, we mention \cite{MZ94} for the stabilization. The first local controllability results on the bilinear Schrödinger equation appear in \cite{Beauchard05,BC06}. These local controllability results have been extended under weaker assumptions in \cite{BL10,BM14}, in a more general setting in infinite time \cite{NN12},  and also in the case of simultaneous controllability of a finite number of particles \cite{Morancey14}, and the references therein. Note that, despite the infinite speed of propagation, it was proved that a minimal amount of time is required for the controllability of some bilinear Schrödinger equations, see \cite{Coron06,BCT14,BM14} (see also \cite{Morancey14,Bournissou23}) and the references therein.  In addition to the exact controllability and the stabilization, various aspects of the controllability of the bilinear Schr\"odinger systems have been investigated. Concerning the approximative controllability, this has been studied by the geometric control techniques via appropriate Galerkin approximations, see e.g., \cite{CMSB09,BCCS12,BCS14}  and the references therein. The Lyapunov technique has been used to obtain the global controllability results, see, e.g., \cite{Mirrahimi09,BM09,Nersesyan09,Nersesyan10} though no indication of the convergence rate is given.

\subsection{The organization of the paper}
The paper is organized as follows. In \Cref{sect-pre}, we establish several results on $A$ and $\cD(A)$. In particular, we prove that $A$ is skew-adjoint in \Cref{lem-A}. \hmb{\Cref{sect-Q} is devoted to various properties of the Gramian operator $Q$. }
In \Cref{sect-WP}, we establish the well-posedness and the stability of various linear and nonlinear Schr\"odinger systems. These results will be used in the proof of the main theorems mentioned above. \Cref{sect-Rapid} is devoted to the rapid stabilization, in particular, we prove \Cref{thm1-S,thm2-S} there. In \Cref{sect-Finite}, we study the finite-time stabilization.  We prove \Cref{thm-FT-LN} using estimates on the cost of controls for the linearized system established there (see \Cref{pro-cost-upper,pro-cost-lower}). The analysis of the upper bound (\Cref{pro-cost-upper}) is based on the moment method. The analysis of the lower bound (\Cref{pro-cost-lower}) is based on a lower bound of the cost of a singular perturbation control problem (\Cref{pro-S-M}). In 
\Cref{sect-appendixA}, we discuss a well-posedness result on control systems associated with operator semi-groups, which is used throughout the paper.

\section{Preliminaries}\label{sect-pre}

In this section, we will prove some properties related to $A$ defined in \eqref{def-A} and $Q$ defined in \eqref{def-Q}. We begin with \hmr{the following result.}

\begin{lemma} \label{lem-density} Let $\gamma_1, \gamma_2 \in \mR$, and let $\cA: \cD(\cA) \subset \mH \to \mH$ be defined by 
\be \label{def-cA}
\cA y = 
\left(\begin{array}{cc} - \Delta y_2 + \gamma_1 y_2 \\[6pt]
\Delta y_1 + \gamma_2 y_1
\end{array}\right)
 \quad \mbox{ with } \quad \cD(\cA) = \Big\{ y \in \mH; \cA y \in \mH \Big\}.  
\ee
We have 
\begin{itemize}
\item[i)] The set $\cD(\cA^\infty)$ is dense in $\mH$. 
\item[ii)] The set $\cD(\cA^\infty) \cap \mH_{1, \sharp}$ is dense in $\mH_{1, \sharp}$. 
\item[iii)] The set $\cD(\cA^\infty) \cap \mH_{1, \sharp}$ is dense in $\cD (\cA) \cap \mH_{1, \sharp}$ equipped \hmr{with} the graph-norm of $\cD(\cA)$. 
\end{itemize}
\end{lemma}

Recall that 
$$
\cD(\cA^\infty) = \bigcap_{k \ge 1} \cD(\cA^k). 
$$

\begin{proof} We first prove $i)$. Let $y = (y_1, y_2)\tr \in \mH$. Then 
\be \label{lem-density-y}
y_1 = \sum_{k = 1}^\infty a_k \varphi_k \quad \mbox{ and } \quad y_2 = \sum_{k = 1}^\infty b_k \varphi_k, 
\ee
for some $(a_k), (b_k) \subset \mR$ such that $\sum_{k \ge 1} \lambda_k^3 (|a_k|^2 + |b_k|^2) < + \infty$.  Denote 
\be \label{lem-density-yn}
y_{1, n} = \sum_{k = 1}^n a_k \varphi_k \quad \mbox{ and } \quad y_{2, n} = \sum_{k = 1}^n b_k \varphi_k,  
\ee
and set 
$$
y_n = (y_{1, n}, y_{2, n})\tr. 
$$

Since 
\be \label{lem-density-varphik}
\varphi_k'' = - \lambda_k \varphi_k \mbox{ in } I \quad \mbox{ and } \quad \varphi_k = 0 \mbox{ on } \partial I, 
\ee
it follows from \eqref{CGM1} that 
\be \label{lem-density-p1-1}
y_n \in \cD(\cA^\infty). 
\ee
It is clear that  
\be \label{lem-density-p2-1}
y_n = (y_{1, n}, y_{2, n})\tr \to y \mbox{ in } \mH.  
\ee
Assertion $i)$ now follows from \eqref{lem-density-p1-1} and \eqref{lem-density-p2-1}. 

\medskip 
We next deal with $ii)$.  We first note that 
$$
\cD(\cA^\infty) \cap \mH_{1, \sharp} = \left\{ y = (y_1, y_2)\tr \in \cD (\cA^\infty);   \int_I y_1 \varphi_1 = 0  \right\}
$$
Let $y = (y_1, y_2)\tr \in \mH_{1, \sharp}$ and define $y_n$ by \eqref{lem-density-yn} using \eqref{lem-density-y}. Then  
\be \label{lem-density-p1}
y_n \in \cD(\cA^\infty) \quad \mbox{ and } \quad y_n \to y \mbox{ in } \mH. 
\ee
Define $\hy_n = (\hy_{1, n}, \hy_{2, n})\tr $ by
\be \label{lem-density-p2}
\hy_{1, n} = y_{1, n} - \langle y_{1, n}, \varphi_1 \rangle_{L^2(I)} \varphi_1 \quad \mbox{ and } \quad \hy_{2, n} = y_{2, n}.  
\ee

It follows from \eqref{lem-density-varphik} that 
\be \label{lem-density-hyn}
\hy_n \in \cD(\cA^\infty) \cap \mH_{1, \sharp}. 
\ee
Using the fact  
\be \label{lem-density-fact}
\langle y_{1, n}, \varphi_1 \rangle_{L^2(I)} \to \langle y_{1}, \varphi_1 \rangle_{L^2(I)} \mathop{=}^{y \in \mH_{1, \sharp}} 0, 
\ee
we derive from \eqref{lem-density-p1} and \eqref{lem-density-p2} that 
$$
\hy_n \to y \mbox{ in } \mH. 
$$
Assertion $ii)$ is proved. 

\medskip 
We finally establish $iii)$.   Let $y = (y_1, y_2)\tr \in \cD(\cA) \cap \mH_{1, \sharp}$ and define $y_n$ by \eqref{lem-density-yn} using \eqref{lem-density-y}. Then  
\be \label{lem-density-ccccc}
y_n \in \cD(\cA^\infty) \quad \mbox{ and } \quad y_n \to y \mbox{ in } \cD(\cA). 
\ee
Define $\hy_n = (\hy_{1, n}, \hy_{2, n})\tr $ by \eqref{lem-density-p2}. Then, by \eqref{lem-density-hyn},  
$$
\hy_n \in \cD(\cA^\infty) \cap \mH_{1, \sharp}. 
$$
Using \eqref{lem-density-fact}, we derive from \eqref{lem-density-ccccc} that 
$$
\hy_n \to y \mbox{ in } \cD(\cA). 
$$
Assertion $iii)$ is established.

\medskip 
The proof is complete.  
\end{proof}

We next establish a result which implies that $A$ is skew-adjoint. 

\begin{lemma} \label{lem-A}  Let $\gamma \in \mR$ and let $\cA: \cD(\cA) \subset \mH \to \mH$ be defined by 
\be \label{def-cA}
\cA y = 
\left(\begin{array}{cc} - \Delta y_2 - \gamma y_2 \\[6pt]
\Delta y_1 + \gamma y_1
\end{array}\right)
 \quad \mbox{ and } \quad \cD(\cA) = \Big\{ y \in \mH; \cA y \in \mH \Big\}.  
\ee
Then $\cA$ is skew-adjoint, i.e., $\cD(\cA^*) = \cD(\cA)$ and $\cA^* = -\cA$ in $\cD(\cA)$. 
\end{lemma}

We recall, by \Cref{lem-density}, that $\cD(\cA)$ is dense in $\mH$. 

\begin{proof} Since 
$$
\varphi_k'' = - \lambda_k \varphi_k \mbox{ in } I, 
$$
we derive from the definition of $\mH$ that $y  = (y_1, y_2)\hmr{\tr} \in \cD(\cA)$ if and only if 
\be\label{lem-A-p1}
y  \in [H^5(I)]^2, 
\ee
\be\label{lem-A-p2}
y_\ell(x) =  y_\ell''(x)  = y_\ell''''(x)  = 0 \mbox{ for } x \in \partial I, \ell = 1, 2.  
\ee
Using this fact, we derive, by integration by parts, for $y = (y_1, y_2) \tr  \in \cD(\cA)$ and $z = (z_1, z_2)\hmr{\tr}  \in \cD(\cA)$, that 
\be\label{lem-A-p1}
- \int_I  \Delta y_2 z_1 = -\int_{I} y_2 \Delta z_1, \quad - \int_I  \Delta y_2' z_1' = -\int_{I} y_2' \Delta z_1', 
\ee
\be\label{lem-A-p2}
- \int_I  \Delta y_2'' z_1'' = -\int_{I} y_2'' \Delta z_1'', \quad - \int_I  \Delta y_2''' z_1''' = -\int_{I} y_2''' \Delta z_1''', 
\ee
\be\label{lem-A-p3}
\int_I  \Delta y_1 z_2 = \int_{I} y_1 \Delta z_2, \quad  \int_I  \Delta y_1' z_2' = \int_{I} y_1' \Delta z_2', 
\ee
\be\label{lem-A-p4}
\int_I  \Delta y_1'' z_2'' = \int_{I} y_1'' \Delta z_2'', \quad  \int_I  \Delta y_1''' z_2''' = \int_{I} y_1''' \Delta z_2'''. 
\ee
It follows that, for $y \in \cD(\cA)$ and $z \in \cD(\cA)$,  
$$
\langle \cA y, z \rangle_{\mH} = \langle y, - \cA z \rangle_{\mH}.
$$

It remains to show that $\cD(A^*) \subset \cD(A)$. This is equivalent to establish 
that  if  $z \in \mH$ is  such that  
\be \label{lem-A-p5}
|\langle \cA y, z \rangle_{\mH} | \le C \| y\|_{\mH} \mbox{ for all } y \in \cD(\cA) 
\ee
for some positive constant $C = C(z)$ independent of $y$, then $z \in \cD(\cA)$. 

Indeed, fix such a $z$. From \eqref{lem-A-p5}, we deduce from \eqref{lem-A-p1} and \eqref{lem-A-p2} that, for $y \in \cD(\cA)$,  
\be \label{lem-A-ineq}
\left| - \int_I  \Delta y_2''' z_1''' + \int_I  \Delta y_1''' z_2''' \right| \le C \| y\|_{\mH}. 
\ee
By taking $y_1 = 0$ in \eqref{lem-A-ineq}, we obtain that, for $y_2 \in H^5(I)$ with $y_2 = y_2'' = y_2''''=0$ on $\partial I$, it holds
\be \label{lem-A-pp1}
\left| - \int_I  \Delta y_2''' z_1''' \right| \le C \| y_2 \|_{H^3(I)}. 
\ee

Given $\varphi \in C^\infty(\bar I; \mR)$, define, for $x \in \bar I$, 
\be \label{lem-A-xi3}
 \xi_3 (x) =  \varphi (x) - \frac{1}{2} \varphi'(1) x^2 + \frac{1}{2} \varphi'(0) (1 - x)^2, 
\ee
$$
\quad \xi_2 (x) = \int_0^x \xi_3(s) \, ds - x \int_0^1 \xi_3(s) \, ds,    \quad 
\xi_1(x) = \int_0^x \xi_2(s) \, ds, 
$$
and 
\be\label{lem-A-y2}
y_2 (x) = \int_0^x \xi_1(s) \, ds - x \int_0^1 \xi_1 (s) \, ds. 
\ee
Simple computations give, for $x \in \bar I$, 
$$
y_2' = \xi_1 - \int_0^1 \xi_1 (s) \, ds, \quad 
y_2'' = \xi_2,  \quad y_2''' = \xi_3 - \int_0^1 \xi_3(s) \, ds, \quad y_2'''' = \varphi' - \varphi'(1) x - \varphi'(0) (1 - x).  
$$

One can then check that $y_2 \in H^5(I)$ with $y_2 = y_2'' = y_2''''=0$ on $\partial I$. It follows from \eqref{lem-A-pp1} applied to $y_2$ given by \eqref{lem-A-y2} that 
$$
\left| - \int_I  \big(\Delta \varphi  - \varphi'(1) + \varphi'(0) \big) z_1''' \right| \le C \| \xi_3 \|_{L^2(I)} \mbox{ for } \varphi \in C^\infty (\bar I; \hmr{\mR}),    
$$
where $\xi_3$ is defined by \eqref{lem-A-xi3}.  Since $z_1'' = 0$ on $\partial I$, we deduce that  
$$
\left| - \int_I  \Delta \varphi z_1''' \right| \le C \| \xi_3 \|_{L^2(I)} \mbox{ for } \varphi \in C^\infty (\bar I; \hmr{\mR}),    
$$
where $\xi_3$ is defined by \eqref{lem-A-xi3}. By first considering $\varphi \in C^\infty_{c} (I)$ and then using $\varphi \in C^\infty(\bar I)$ with $\varphi' = 0$ on $\partial I$, we obtain 
\be\label{lem-A-z1}
z_1''' \in H^2(I) \quad \mbox{ and } \quad z_1'''' = 0 \mbox{ on } \partial I. 
\ee

Similarly, by taking $y_2 = 0$ in \eqref{lem-A-ineq}, we derive that 
\be \label{lem-A-z2}
z_2''' \in H^2(I) \quad \mbox{ and } \quad z_2'''' = 0 \mbox{ on } \partial I. 
\ee

Combining \eqref{lem-A-z1} and \eqref{lem-A-z2} yields 
$$
z \in \cD(\cA). 
$$

The proof is complete. 
\end{proof}

\section{Properties of the Gramian operator} \label{sect-Q}

In this section, we study properties of $Q = Q (\lambda)$. Recall that $A$ is defined in \eqref{def-A}, $B^*$ is given by \eqref{def-B*},  and $\varphi_k$ is given in \eqref{def-lambdak-varphik}.  Set, for $k \ge 1$,  
$$
\Phi^1_k =  \left( \begin{array}{c}
\varphi_k \\[6pt]
0
\end{array} \right) \quad \mbox{ and } \quad \Phi^2_k =  \left( \begin{array}{c}
0 \\[6pt]
\varphi_k
\end{array} \right). 
$$
Thus $\Phi^1_1 = \Phi_1$.  We begin with the following elementary fact which is a direct consequence of the definition of $A$ and $\varphi_k$. 

\begin{lemma} \label{lem-Q2} 
We have, for $k \ge 1$ and $t \ge 0$,  
\be \label{formula-St}
e^{t A} \Phi^1_k =  \varphi_k \left( \begin{array}{c}
\cos [(\lambda_k - \lambda_1) t] \\[6pt]
- \sin [(\lambda_k - \lambda_1) t]
\end{array} \right) \quad \mbox{ and } \quad e^{t \cA} \Phi^2_k = \varphi_k \left( \begin{array}{c}
\sin [(\lambda_k - \lambda_1) t] \\[6pt]
\cos [(\lambda_k - \lambda_1) t]
\end{array} \right). 
\ee
\end{lemma}

Using \Cref{lem-Q2}, we can prove the following result. 

\begin{lemma} \label{lem-Q3} Let $\lambda > 0$. For $i, j = 1, 2$, and for $k, l \ge 1$, set 
$$
d^{i, j}_{k, l} := \langle Q  \Phi^i_k,   \Phi^j_l \rangle_{\mH} =  \int_{0}^\infty e^{- 2 \lambda s} \langle B^*e^{s \cA} \Phi^i_k, B^*e^{s \cA} \Phi^j_l \rangle_{\mR} \, ds. 
$$
We have, for some positive constant $\Lambda \ge 1$,  
\be \label{lem-Q3-cl0-1}
 \quad d^{1, j}_{1, l} = 0 \mbox{ for } j =1, 2, l \ge 1, 
\ee
\be \label{lem-Q3-cl0-2}
|d^{i, j}_{k, k}| \le \frac{\Lambda k^6}{\lambda_k} \mbox{ for } i, j =1, 2, i \neq j, k \ge 1, \quad 
|d^{i, j}_{k, l}| \le  \frac{\Lambda k^3 l^3}{|\lambda_k - \lambda_l|} \mbox{ for } i, j = 1, 2, k, l \ge 1,  k \neq l.  
\ee
\be \label{lem-Q3-cl0-3}
|d^{i, i}_{k, k}| \le \Lambda k^6 \mbox{ for } i = 1, 2, k \ge 2, \quad \mbox{ and }  \quad |d^{2, 2}_{1, 1}| \le \Lambda.  
\ee
As a consequence, 
\be\label{lem-Q3-cl1}
Q: \mH_{1, \sharp} \to \mH_{1, \sharp}. 
\ee
Assume in addition that \eqref{cond-mu-EEE} holds. We have
\be \label{lem-Q3-cl2}
 \Lambda^{-1} k^6 \le |d^{i, i}_{k, k}| \le \Lambda k^6 \mbox{ for } i = 1, 2, k \ge 2, \quad \mbox{ and } \quad  \Lambda^{-1} \le |d^{2, 2}_{1, 1}| \le \Lambda.  
\ee
\end{lemma}

\begin{proof} Using \eqref{formula-St}, we derive from \eqref{def-B*} that, for $k \ge 1$ and $s \ge 0$,  
\be \label{lem-Q3-p1}
B^*e^{s \cA} \Phi^1_k =  - \sin[(\lambda_k - \lambda_1) s] \Big( \langle \mu \varphi_1,  \varphi_k \rangle_{H^3(I)}  - (\mu \varphi_1)_{xx} (1)  \varphi_{k, xxx}(1) 
+ (\mu \varphi_1)_{xx}(0)  \varphi_{k, xxx}(0) \Big), 
\ee
and 
\be \label{lem-Q3-p2}
B^*e^{s \cA} \Phi^2_k =  \cos[(\lambda_k - \lambda_1) s] \Big( \langle \mu \varphi_1,  \varphi_k \rangle_{H^3(I)}  - (\mu \varphi_1)_{xx} (1)  \varphi_{k, xxx}(1) 
+ (\mu \varphi_1)_{xx}(0)  \varphi_{k, xxx}(0) \Big). 
\ee
Note that, for some positive constant $C \ge 1$,  
\be \label{lem-Q3-p3}
\Big| \Big( \langle \mu \varphi_1,  \varphi_k \rangle_{H^3(I)}  - (\mu \varphi_1)_{xx} (1)  \varphi_{k, xxx}(1) 
+ (\mu \varphi_1)_{xx}(0)  \varphi_{k, xxx}(0) \Big) \Big| \le C k^3. 
\ee
Assertions \eqref{lem-Q3-cl0-1}, \eqref{lem-Q3-cl0-2}, and \eqref{lem-Q3-cl0-3} now follow from  \eqref{lem-Q3-p1}, \eqref{lem-Q3-p2}, and \eqref{lem-Q3-p3}. 

To derive \eqref{lem-Q3-cl2}, one just needs to note that  \eqref{cond-mu-EEE-cl} holds by \eqref{cond-mu-EEE} (see \Cref{lem-Obs}). The proof is complete. 
\end{proof}

We are ready to introduce the linear application 
$$
T: \mH_{1, \sharp} \to \mH_{1, \sharp} 
$$
defined by 
$$
T (\Phi^i_k) = \frac{d^{i,i}_{k, k}}{\|\Phi^i_k \|^2_{H^3(I)}}  \Phi^i_k \mbox{ for } (i =1, 2 \mbox{ and } k \ge 2) \mbox{ and for } (i=2 \mbox{ and } k = 1).  
$$

Using \Cref{lem-Q3}, we obtain the following result on $T$, which is a direct consequence of \eqref{lem-Q3-cl2}.  
 
\begin{lemma} \label{lem-Q4} Assume \eqref{cond-mu-EEE}. We have, for $s \ge 0$,  
$$
T:  \mH_{1, \sharp}  \cap \big[ H^{3 + s}(I) \big]^2 \to \mH_{1, \sharp}  \cap \big[ H^{3 + s}(I) \big]^2
$$
is linear and invertible. 
\end{lemma}


We also have the following result. 

\begin{lemma} \label{lem-Q5} Let $1/2  <   s < s' \le 1$.  We have 
$$
Q - T:  \mH_{1, \sharp} \cap \big[ H^{3 + s}(I) \big]^2  \to \mH_{1, \sharp}  \cap \big[ H^{3 + s'}(I) \big]^2
$$
is linear and continuous. 
\end{lemma}

\begin{proof} Let $\Phi \in \mH_{1, \sharp} \cap \big[ H^{3 + s}(I) \big]^2$. We represent $\Phi$ under the form 
$$
\Phi = \sum_{j=1, 2, k \ge 1} a^{j, k} \Phi^j_k, 
$$
for some $a^{j, k} \in \mR$ with $a^{1, 1} = 0$ and 
$$
\sum_{j=1, 2, k \ge 1} k^{2 (3 + s)}|a^{j, k}|^2 < + \infty. 
$$

Set 
$$
\beta_l = \langle \varphi_l, \varphi_l \rangle_{\mH} \mbox{ for } l \ge 1. 
$$
Then 
\be \label{lem-Q4-p1}
\beta_l \ge C l^6 \mbox{ for } l \ge 1. 
\ee
Since $\Big(\frac{1}{\beta_k^{1/2}}\Phi^{i}_k \Big)_{i=1, 2, k \ge 1}$ is an orthogonal of $\mH$, we have, by \Cref{lem-Q3}, 
\begin{multline*}
(Q-T)\left(\sum_{j=1, 2, k \ge 1} a^{j, k} \Phi^j_k \right) = \sum_{j=1, 2, k \ge 1}  a^{j, k} \sum_{l \ge 1, l \neq k, i=1, 2} \frac{1}{\beta_l} d^{j, i}_{k, l} \Phi^i_l + \sum_{j=1, 2, k \ge 1} \sum_{i=1,2, i \neq j} a^{j, k} \frac{1}{\beta_k} d^{j,i}_{k, k} \Phi^i_k \\[6pt]
= \sum_{l \ge 1, i =1, 2} \sum_{j=1, 2, k \ge 1, k \neq l}  \frac{a^{j, k} }{\beta_l} d^{j, i}_{k, l} \Phi^i_l + \sum_{j=1, 2, l \ge 1} \sum_{i=1,2, i \neq j} a^{j, l} \frac{1}{\beta_l} d^{j,i}_{l, l} \Phi^i_l. 
\end{multline*}
This implies 
\begin{multline*}
\left\| (Q-T)\Big(\sum_{j=1, 2, k \ge 1} a^{j, k} \Phi^j_k \Big) \right\|_{H^{3+s'}(I)}^2 \mathop{\le}^{\Cref{lem-Q2}} C \sum_{l \ge 1, i =1, 2} \left( \sum_{j=1, 2, k \ge 1, k \neq l}  \frac{|a^{j, k}|}{\beta_l} |d^{j, i}_{k, l}| l^{3+s'} \right)^2 \\[6pt] + 
C \sum_{j=1, 2, l \ge 1} |a^{j, l}|^2 l^{2(1 + s')}
 \\[6pt]
\mathop{\le}^{\Cref{lem-Q2}, \eqref{lem-Q4-p1}} C \sum_{l \ge 1, i =1, 2} \left( \sum_{j=1, 2, k \ge 1, k \neq l}  \frac{|a^{j, k}| k^3 l^{s'}}{|\lambda_k- \lambda_l|} \right)^2 + C \sum_{j=1, 2, l \ge 1} |a^{j, l}|^2 l^{2(1 + s')} \\[6pt] 
\le C \sum_{l \ge 1, i =1, 2} \left(  \sum_{j=1, 2, k \ge 1, k \neq l} |a^{j, k}|^2 k^{2(3 + s)} \right) \left(  \sum_{k \neq l} \frac{l^{2s'}}{k^{2s} (\lambda_k - \lambda_l)^2} \right) + C \sum_{j=1, 2, l \ge 1} |a^{j, l}|^2 l^{2(1 + s')}. 
\end{multline*}
We thus obtain
\be\label{lem-Q4-p2}
\left\| (Q-T) ( \Phi ) \right\|_{H^{3+s'}(I)}^2  \le C \left\| \Phi \right\|_{H^{3+s}(I)}^2 \left( 1 +  \sum_{l \ge 1}\sum_{k \neq l} \frac{l^{2s'}}{k^{2s} (\lambda_k - \lambda_l)^2} \right). 
\ee

Since, for $k \neq l$,  
$$
(\lambda_k - \lambda_l)^2 \le C (k-l)^2 (k+ l^2), 
$$
$$
\sum_{k \neq l} = \sum_{k \le l /2} + \sum_{l/2 < k \le 2 l, k \neq l} + \sum_{k >  2 l}, 
$$
and $1/2 <  s  < s' \le 1$, 
it follows that 
\begin{multline}\label{lem-Q4-p3}
\sum_{l \ge 1} \sum_{k \neq l} \frac{l^{2s'}}{k^{2s} (\lambda_k - \lambda_l)^2} 
\le C\sum_{l \ge 1} \sum_{k \neq l} \frac{l^{2s'}}{k^{2s} (k - l)^2 (k+l)^2} \\[6pt]
  \le C \sum_{l \ge 1} (l^{2s' - 4} + l^{2s' - 2s - 2} + l^{2s' - 3 - 2 s} ) \le C.  
\end{multline}
Combining \eqref{lem-Q4-p2} and \eqref{lem-Q4-p3} yields 
$$
\|(Q - T)(\Phi) \|_{H^{3 + s'}(I)} \le C \| \Phi \|_{H^{3 + s}(I)}. 
$$
The proof is complete. 
\end{proof}

As a consequence of \Cref{lem-Q4} and \Cref{lem-Q5}, and the Fredholm theory, we have the following result. 
\begin{lemma} \label{lem-Q6} Let $1/2 < s< 1$. We have
$$
Q: \mH_{1, \sharp} \cap  \big[ H^{3 + s}(I) \big]^2  \to \mH_{1, \sharp}  \cap  \big[ H^{3 + s}(I) \big]^2
$$
is invertible. We also have 
$$
\| \cQ^{-1}(\Phi^{i}_k) \|_{C^3(\bar I)} \le C k^3 \mbox{ for } (i \ge 2, k \ge 1) \mbox{ and for } (i=2 \mbox{ and } k=1).   
$$
\end{lemma}

\begin{proof}  The first part is a consequence of the Fredholm theory and the three facts:

$a)$ $T + (Q-T): \mH_{1, \sharp} \cap \big[ H^{3 + s}(I) \big]^2 \to \mH_{1, \sharp}  \cap \big[ H^{3 + s}(I) \big]^2 $ is injective since $Q: \mH_{1, \sharp} \to \mH_{1, \sharp} $ is bijective; 

$b)$  $T: \mH_{1, \sharp} \cap \big[ H^{3 + s}(I) \big]^2  \to \mH_{1, \sharp}  \cap \big[H^{3 + s}(I) \big]^2$ is invertible by \Cref{lem-Q4};

$c)$ $Q-T: \mH_{1, \sharp} \cap \big[ H^{3 + s} (I) \big]^2 \to \mH_{1, \sharp}  \cap \big[ H^{3 + s}(I) \big]^2$ is compact by \Cref{lem-Q5}.

\medskip 
We next prove the second part. Note that 
$$
Q \left(\frac{d^{i,i}_{k, k}}{\|\Phi^i_k \|^2_{H^3(I)}}\cQ^{-1}(\Phi^{i}_k) - \Phi^{i}_k \right) = (T - Q)\Phi^{i}_k. 
$$
We have, for $1/2 < s < 1$,  
\be
\|(T - Q)(\Phi^{i}_k) \|_{H^{3 + s}(I)}^2 \mathop{\le}^{\Cref{lem-Q2}}  C \sum_{l \ge 1, l \neq k} \frac{k^{6} l^{12+2s}}{l^{12} (k-l)^2 (k+l)^2} + C k^{2(1 + s)}. 
\ee
Since
$$
\sum_{l \neq k} = \sum_{l \le k/2} + \sum_{k/2 < l \le 2 k, l \neq k} + \sum_{l >  2 k}, 
$$
it follows that, for $1/2 < s < 1$,  
$$
\|(T - Q)\Phi^{i}_k \|_{H^{3 + s}(I)}^2 \le  C \Big( k^{2 + 2 s + 1} + k^{4+ 2s} + k^{6+2s - 3} \Big)+ C k^{2(1 + s)} \le C k^{4 + 2s}. 
$$
Applying the first part, we derive that, for $1/2 < s < 1$,   
$$
\left\|\frac{d^{i,i}_{k, k}}{\|\Phi^i_k \|^2_{H^3(I)}}\cQ^{-1}(\Phi^{i}_k) - \Phi^{i}_k \right\|_{H^{3 + s}(I)} \le Ck^{2 + s}.  
$$
This implies, by the Sobolev embedding $H^{s}(I)$ into $C(\bar I)$,  
$$
\left\| \frac{d^{i,i}_{k, k}}{\|\Phi^i_k \|^2_{H^3(I)}}\cQ^{-1}(\Phi^{i}_k) - \Phi^{i}_k  \right\|_{C^{3}(\bar I)} \le C k^{2 + s}, 
$$
and the conclusion follows. 
\end{proof}

\section{Well-posedness and stability of Schr\"odinger systems} \label{sect-WP}

In this section, we establish the well-posedness and the stability of various systems related to the linear system \eqref{sys-LN2} and the nonlinear system \eqref{sys-NL2}. The main goal is to formulate and establish results which are compatible with the theory of control systems associated with semi-group. Without the language of semi-group, some related results can be found in \cite{BL10}. 

We first introduce $ \bA : \cD(\bA) \subset \bH \to \bH $ defined by 
\be
\bA \Psi = i \Delta \Psi \quad \mbox{ and } \quad \cD(\bA) = \Big\{\Psi \in \bH; \bA \Psi \in \bH \Big\}. 
\ee

\hmr{The following properties holds for the operator ${\bA}$. }

\begin{lemma}\label{lem-bA} We have 
$$
\mbox{$\cD(\bA)$ is dense in $\bH$ and $\bA$ is skew-adjoint.}
$$ 
\end{lemma}

\begin{proof} \rm The conclusion is a consequence of \Cref{lem-density,lem-A} with $\gamma = 0$ after considering the real part and the imaginary part of $\Psi$ and $\bA \Psi$. 
\end{proof}

We next introduce a useful operator related to the definitions of $B$ in \eqref{def-B} and $u F(y)$ in \eqref{def-uF}. 

\begin{definition} \label{def-T} Given $T>0$.  Define 
$$
\bT: L^2( (0, T); H^3(I; \mC) \cap H^1_0 (I; \mC)) \to L^1((0, T); \cD(\bA^*)')
$$
by, for all $\varphi \in \cD(\bA^*)$,  
\be\label{lem-WP-f1}
\langle \bT(f)(t, \cdot), \varphi \rangle_{\cD(\bA^*)', \cD(\bA^*)} = \langle f(t, \cdot), \varphi \rangle_{H^3(I)}  -  f_{xx}(t, 1) \hmr{\overline{\varphi_{xxx} (1)}} + f_{xx}(t, 0) \hmr{\overline{\varphi_{xxx} (0)}},   
\ee
for $f \in L^2( (0, T); H^3(I; \mC) \cap H^1_0 (I; \mC))$. 

\end{definition}

\hmr{It follows from the definition of  $\bH$ that $\varphi \in H^5(I)$ if $\varphi \in \cD(\bA)$ and $\| \varphi \|_{H^5(I)} \le C \| \varphi \|_{\cD(\bA)}$. Thus $\bT$ is  well-defined. The motivation of the definition of $\bT$ will be clear from \Cref{lem-product} below.}

\medskip
We next discuss the well-posedness and the stability of linear systems. \hmb{It is convenient to introduce the space $\bH_{rel}$ which is defined by 
$$
\bH_{rel} = \big\{y \in \bH; y \mbox{ is real} \big\}. 
$$
We also equip this space with the scalar product from $\bH$. Then $\bH_{rel}$ is a real Hilbert space.  We have the following important result on the well-posedness of the linearized system and the regularizing effect. }

\begin{lemma} \label{lem-WP} Let $0< T < T_0$ and  $\lambda \in \mR$. Let $\Phi_0 \in \bH$ and  $\hmr{\ff} \in L^1((0, T); \cD(\bA^*)')$.  There exists a unique weak solution $\Phi \in C([0, T]; \cD(\bA^*)\hmr{'})$ to the system 
\be\label{lem-WP-sys}
\left\{\begin{array}{cl}
i \Phi_t = - \Delta \Phi  - \lambda \Phi   + \hmr{\ff}  & \mbox{ in } (0, T)  \times I, \\[6pt] 
\Phi(t, 0) = \Phi(t, 1) = 0 &  \mbox{ in } (0, T), \\[6pt]
\Phi(0, \cdot) = \Phi_0 & \mbox{ in } I,  
\end{array} \right. 
\ee 
i.e. 
\be \label{def-WS}
i \frac{d}{dt} \langle \Phi, \Psi \rangle_{\bH} = -  \langle \Phi, \Delta \Psi \rangle_{\bH} - \lambda   \langle \Phi, \Psi \rangle_{\bH} + \langle \hmr{\ff}, \Psi \rangle_{\cD(\bA^*)', \cD(\bA^*)} \mbox{ in } (0, T)
\ee
in the distributional sense for all $\Psi \in \cD(\bA^*)$. Let $f \in L^2((0, T); H^3(I; \hmr{\mC}) \cap H^1_0(I; \hmr{\mC}))$. Define $\hmr{\ff} = \bT(f) \in L^1((0, T); \cD(\bA^*)')$. Then the weak solution $\Psi$ of \eqref{lem-WP-sys} satisfies $\Psi \in C([0, T]; \bH)$ and 
\be \label{lem-WP-est}
\|\Phi(t, \cdot)\|_{H^3(I)} \le C\Big( \| \Phi_0\|_{H^3(I)}  + \| f\|_{L^2((0, t); H^3(I))} \Big) 
\mbox{ in } [0, T], 
\ee
for some positive constant $C$ depending only on $T_0$ and $\alpha$.   \hmb{Let $\bR_{ij} \in \cL (\bH_{rel})$ for $i, j = 1, 2$
be such that $\|\bR_{ij} \varphi_k\|_{C^3(\bar I)} \le \alpha k^3$ for some $\alpha > 0$ and for all $k \ge 1$. Denote 
$$
\bR \varphi = (\bR_{11} \varphi_1 + \bR_{12} \varphi_2) + i (\bR_{21} \varphi_1 + \bR_{22} \varphi_2) \mbox{ for } \varphi \in \bH, 
$$
where $\varphi_1 = \Re \varphi $ and $\varphi_2 = \Im \varphi$. We have, $(\bR \Phi)_{xxx} \in C([0, 1]; H^3(0, T))$, and 
\be \label{lem-WP-reg}
\| (\bR \Phi)_{xxx} (\cdot, x) \|_{L^2(0, T)} \le C_1 \| \Phi_0\|_{H^3(I)}  + C_2 \| f\|_{L^2((0, T),  H^3(I))},  
\ee
for some positive constants $C_1, C_2$ depending only on $T_0$ and $\alpha$. 
Moreover, $C_2$ can be chosen such that $C_2 \to 0$ as $T_0 \to 0_+$.}
\end{lemma}

\begin{proof}[Proof of \Cref{lem-WP}] The existence and uniqueness of solutions in $C([0, T]; \cD(\bA^*)')$ follows from \Cref{pro-WP} in the appendix. \hmr{We next} show the existence of a solution $\Psi \in C([0, T]; \bH)$ satisfying \eqref{lem-WP-est} \hmr{in the case $\ff = \bT (f)$}. 

\medskip
We first deal with the system 
\be\label{lem-WP-sys-aux}
\left\{\begin{array}{cl}
i \Phi_t = - \Delta \Phi + \hmr{\ff}  & \mbox{ in } (0, T) \times I, \\[6pt] 
\Phi(t, 0) = \Phi(t, 1) = 0 &  \mbox{ in } (0, T), \\[6pt]
\Phi(0, \cdot) = \Phi_0 & \mbox{ in } I
\end{array} \right. 
\ee 
instead of \eqref{lem-WP-sys}, i.e., we consider \eqref{lem-WP-sys} with $\lambda = 0$.  We search $\Phi \in C([0, T]; \bH)$ under the form 
$$
\Phi(t, x) = \sum_{k \ge 1} a_k (t) \varphi_k (x) \mbox{ in } (0, T) \times I. 
$$
Using \hmr{\eqref{def-WS}} with $\lambda =0$ and $\varphi= \varphi_k$, we obtain 
\be \label{lem-WP-ak0}
i (1 + \lambda_k + \lambda_k^2 + \lambda_k^3) a_k' = \lambda_k (1 + \lambda_k + \lambda_k^2 + \lambda_k^3) a_k + c_k \mbox{ in } (0, T),   
\ee
where 
\be \label{lem-WP-bk}
c_k(t) =  \langle f(t, \cdot),  \varphi_{k} \rangle_{H^3(I)}  -  f_{xx}(t, 1) \varphi_{k, xxx} (1) + f_{xx}(t, 0) \varphi_{k, xxx} (0) \mbox{ in } (0, T). 
\ee
We derive from \eqref{lem-WP-ak0} that 
$$
a_k' = - i \lambda_k a_k - i b_k \mbox{ in } (0, T) \mbox{ where } b_k = \frac{c_k}{1 + \lambda_k + \lambda_k^2 + \lambda_k^3}. 
$$
We then get 
\be \label{lem-WP-ak}
a_k (t) = e^{- i \lambda_k t} a_k(0) -   i \int_0^t e^{- i \lambda_k (t-s)} b_k (s) \mbox{ in } (0, T).     
\ee

Combining \eqref{lem-WP-bk} and \eqref{lem-WP-ak} yields 
\begin{multline}\label{lem-WP-abk}
\sum_{k \ge 1} k^6 |a_k (t)|^2 \le C  \| \Phi_0\|_{H^3(I)}^2 + 
C  \sum_{k \ge 1} \frac{1}{\lambda_k^3}\int_0^t |\langle f, \varphi_k \rangle_{H^3(I)}|^2 \, dt \\[6pt]
+ C \sum_{k \ge 1} \Big| \int_0^t e^{i \lambda_k s} f_{xx}(s, 1) \, ds \Big|^2 
+ C \sum_{k \ge 1} \Big| \int_0^t e^{i \lambda_k s} f_{xx}(s, 0)  \, ds \Big|^2. 
\end{multline}
Here and in what follows in this proof, $C$ denotes a positive constant depending only on $T_0$. We have 
\be \label{lem-WP-ccc1}
\sum_{k \ge 1}  \frac{1}{\lambda_k^3} \int_0^t |\langle f, \varphi_k \rangle_{H^3(I)}|^2 \, dt  \le C \int_0^t \|f(s, \cdot) \|_{L^2((0, T); H^3(I))}^2 \, ds.    
\ee
Applying Ingham's inequality (see, e.g., \cite[Theorem 4.3 on page 59]{KL05}) and using the properties of Riesz basis, see, e.g., \cite[Theorem 9 on page 32]{Young80}, we obtain 
\be \label{lem-WP-ccc2}
 \sum_{k \ge 1} \Big| \int_0^t e^{i \lambda_k s} f_{xx}(s, 1) \, ds \Big|^2 \le C \int_0^t |f_{xx}(s, 1)|^2 \, ds 
\ee
and 
\be \label{lem-WP-ccc3}
\sum_{k \ge 1} \Big| \int_0^t e^{i \lambda_k s} f_{xx}(s, 0)  \, ds \Big|^2 \le C \int_0^t |f_{xx}(s, 0)|^2 \, ds.  
\ee

Combining \eqref{lem-WP-abk}, \eqref{lem-WP-ccc1}, \eqref{lem-WP-ccc2}, and \eqref{lem-WP-ccc3} yields  
$$
\| \Phi(t)\|_{H^3(I)}^2 \le C\left( \| \Phi_0\|_{H^3(I)}^2  +  \int_0^t  \| f(s, \cdot)\|^2_{H^3(I)} \, ds \right)  \mbox{ in } [0, T].   
$$
One can also check that $\Phi \in C([0, T]; \bH)$ is also a weak solution of \eqref{lem-WP-sys}. 
The conclusion in the case $\lambda = 0$ follows.

\medskip 

We next prove \eqref{lem-WP-reg} still for $\lambda = 0$. For notational ease, we only consider the case where $\bR_{12}, \bR_{21}, \bR_{22}$ are 0.  By \eqref{lem-WP-ak}, we have 
$$
(\bR\Phi)_{xxx}(t, x) = \sum_{k \ge 1}  \lambda_k^{3/2} \xi_{k}(x) \Re a_{k}(t) \mbox{ in } [0, T] \times [0, 1], 
$$
where 
\be \label{lem-WP-xik}
\xi_k(x) = \frac{1}{\lambda_k^{3/2}} (\bR_{11} \varphi_{k})_{xxx} (x), \mbox{ which yields }  |\xi_k(x)| \le  \alpha \mbox{ in } [0, 1]. 
\ee
Extend $f$ by $0$ for $(t, x) \in (\mR \setminus [0, T]) \times (0, 1)$ and we still denote this extension by $f$.   We then have
$$
(\bR \Phi)_{xxx}(t, x) = \sum_{k \ge 1} \xi_k(x) \Re \Big(\alpha_{k, 0}(t) +  \alpha_{k, 1}(t)  + \alpha_{k, 2}(t) + \alpha_{k, 3}(t) \Big) \mbox{ in } [0, T] \times [0, 1],   
$$
where 
$$
\alpha_{k, 0}(t) = \lambda_k^{3/2} a_k(0) e^{- i \lambda_k t} \mbox{ in } [0, T], 
$$
$$
\alpha_{k, j}(t) =  \int_0^t e^{- i \lambda_k  (t -s)} \beta_{k, j}(s) \, ds \mbox{ in } [0, T] \mbox{ for } j=1, 2, 3, 
$$
with 
$$
\beta_{k, 1}(t) =  - i \frac{\lambda_k^{3/2} }{1 + \lambda_k + \lambda_k^2 + \lambda_k^3}  \langle f(t, \cdot),  \varphi_k \rangle_{H^3(I)} \mbox{ in } \mR, 
$$
$$
\beta_{k, 2}(t) = - i \frac{\lambda_k^{3} }{1 + \lambda_k + \lambda_k^2 + \lambda_k^3}   f_{xx}(t, 0),  \quad \beta_{k, 3}(t) = - i \frac{\lambda_k^{3} }{1 + \lambda_k + \lambda_k^2 + \lambda_k^3}  f_{xx}(t, 1) \mbox{ in } \mR. 
$$

Set 
\be \label{lem-WP-rho}
\rho_k(t) = e^{-i \lambda_k t} \mathds{1}_{[0, T]} (t) \mbox{ in } \mR, 
\ee
and denote 
\be \label{lem-WP-ta}
\ta_{k, j}= \rho_k *\beta_{k, j}  \mbox{ in } \mR.
\ee
Then 
$$
\ta_{k, j}  = \alpha_{k, j}  \mbox{ in } [0, T]. 
$$

Given an appropriate function $g$ defined in $\mR$, we denote its Fourier transform by $\hat g$. We have, by \eqref{lem-WP-ta},  
$$
\hat \ta_{k, j}(\xi) = \hat \rho_k(\xi) \hat \beta_{k, j} (\xi) \mbox{ in } \mR.  
$$
Since, by \eqref{lem-WP-rho}, 
$$
\hat \rho_k (\xi) = \frac{1}{-i(\xi + \lambda_k)} \Big(e^{- i (\lambda_k + \xi) T} - 1 \Big) \mbox{ in } \mR, 
$$
it follows  that 
\be \label{lem-WP-reg-p1}
|\hat \rho_k (\xi)| \le C  \quad \mbox{ and } \quad |\hat \rho_k (\xi)|  \le \frac{C}{|\xi + \rho_k|} \mbox{ for }  \xi \in \mR. 
\ee

We claim that, for some positive constant $C$, depending only on $T_0$,  
\be \label{lem-WP-reg-claim}
\sum_{k \ge 1} |\hat \rho_k (\xi)| \le C \mbox{ for } \xi \in \mR. 
\ee

Indeed, from \eqref{lem-WP-reg-p1}, we have, for $-\lambda_{k_0+1} \le  \xi < - \lambda_{k_0}$ and $k_0 \ge 2$, 
\begin{multline} \label{lem-WP-reg-p2}
\sum_{k \ge 1} |\hat \rho_k (\xi)| \le C +  C \sum_{1 \le k \le k_0 - 1 } \frac{1}{\lambda_{k_0} - \lambda_k}  + C \sum_{2 k_0 \ge k \ge k_0 + 1 } \frac{1}{\lambda_{k} - \lambda_{k_0}} + C \sum_{k \ge 2 k_0} \frac{1}{\lambda_{k} - \lambda_{k_0}} \\[6pt] 
\le C + C \sum_{1 \le \ell \le k_0 } \frac{1}{k_0 \ell} + C \sum_{k \ge 2 k_0} \frac{1}{\lambda_k} \le C. 
\end{multline}
It is clear from \eqref{lem-WP-reg-p1} that 
\be  \label{lem-WP-reg-p3}
\sum_{k \ge 1} |\hat \rho_k (\xi)| \le C \mbox{ for } \xi \ge - \lambda_2. 
\ee
Claim \eqref{lem-WP-reg-claim} now follows from \eqref{lem-WP-reg-p2} and \eqref{lem-WP-reg-p3}.

Using \eqref{lem-WP-reg-claim}, we obtain 
$$
\sum_{k \ge 1} |\hat \ta_{k, j}(\xi)| \le  C \Big( |\hat f_{xx}(\xi, 0)| + |\hat f_{xx}(\xi, 1)| \Big) \mbox{ for } j=2, 3.  
$$
This implies, by \eqref{lem-WP-xik},   
\be \label{lem-WP-reg-pp1}
\| \sum_{k \ge 1} \xi_k(x) \Re \alpha_{k, j} \|_{L^2 (0, T)} \le C_{\alpha, T_0} \|f \|_{L^2((0, T); H^3(I))} \mbox{ for } j = 2, 3. 
\ee

We also have 
$$
\left(\sum_{k \ge 1} |\hat \rho_k (\xi)| |\hat \beta_{1, k} (\xi)| \right)^2  \le  \sum_{k \ge 1} |\hat \rho_k (\xi)|^2  \sum_{k \ge 1} |\hat \beta_{1, k} (\xi)|^2 \le C  \sum_{k \ge 1} |\hat \beta_{1, k} (\xi)|^2.  
$$
This yields, by \eqref{lem-WP-xik},  
\be \label{lem-WP-reg-pp2}
\| \sum_{k \ge 1} \xi_k(x) \Re \alpha_{k, 1} \|_{L^2 (0, T)} \le C_{\alpha, T_0} \|f \|_{L^2((0, T); H^3(I))}. 
\ee

We claim that, for $(c_k)_{k \ge 1} \subset \mC$ with  $\sum_{k \ge 1} |c_k|^2 < + \infty$, 
\be \label{lem-WP-claim1}
\| \sum_{k \ge 1} c_k e^{-i\lambda_k t} \|_{L^2(0, T)} \le C \left(\sum_{k \ge 1} |c_k|^2 \right)^{1/2}. 
\ee
Indeed, for $\phi \in L^2(0, T)$, 
\be \label{lem-WP-reg-p4}
\sum_{k \ge 1}  \int_0^T c_k e^{-i\lambda_k t} \overline {\phi(t)} \, dt = \sum_{k \ge 1} c_k \overline{\gamma_k}, \mbox{ where }
\gamma_k = \int_0^T \phi(t) e^{i \lambda_k t } \, dt. 
\ee
Since 
$$
\left| \sum_{k \ge 1} c_k \overline{\gamma_k} \right| \le \left( \sum_{k \ge 1} |c_k|^2 \right)^{1/2} \left(\sum_{k \ge 1} |\gamma_k|^2 \right)^{1/2}, 
$$
and, by Ingham's inequality (see, e.g., \cite[Theorem 4.3 on page 59]{KL05}),  
$$
\left(\sum_{k \ge 1} |\gamma_k|^2 \right)^{1/2} \le C \| \phi\|_{L^2(0, T)}, 
$$
we obtain  Claim \eqref{lem-WP-claim1}. 

As a consequence of \eqref{lem-WP-claim1}, we have,  for $(c_k)_{k \ge 1} \subset \mC$ with  $\sum_{k \ge 1} |c_k|^2 < + \infty$, 
\be \label{lem-WP-claim2}
\| \sum_{k \ge 1} c_k e^{i\lambda_k t} \|_{L^2(0, T)} \le C \left(\sum_{k \ge 1} |c_k|^2 \right)^{1/2}. 
\ee

Combining \eqref{lem-WP-claim1} and \eqref{lem-WP-claim2} yields 
\be \label{lem-WP-reg-pp2}
\| \sum_{k \ge 1} \xi_k(x) \Re \alpha_{k, 1} \|_{L^2 (0, T)} \le C_{\alpha, T_0} \|f \|_{L^2((0, T); H^3(I))}. 
\ee

To obtain the behaviour of $C_{2}$ as $T_0$ goes to 0, one just notes that 
$$
|\hat \rho_k (\xi)| \le  4 T_0 \mbox{ if } |\xi + \lambda_k| \le \frac{1}{T_0}
$$ 
and 
$$
|\hat \rho_k (\xi)| \le \frac{2}{|\xi + \lambda_k|} \mbox{ if } |\xi + \lambda_k| \mbox{ for } \xi \in \mR. 
$$

To obtain the conclusion for a general $\lambda$, one first notes that if $\Phi$ is a solution of \eqref{lem-WP-sys-aux} then $\Phi(t, x) e^{- i\lambda t}$ is a solution of \eqref{lem-WP-sys} with the same initial condition and with 
the source $ e^{-i\lambda t} f(t, x)$ and then apply the result in the case $\lambda =0$ to reach the conclusion.  

\medskip 
The proof is complete. 
\end{proof}

The following simple result is useful to compare with previous results and motivates the definition of the operator $\bT$ in \eqref{lem-WP-f1}. 

\begin{lemma} \label{lem-product} Let $0< T < T_0$ and  $\lambda \in \mR$. Let $\Phi_0 \in \bH$ and let $f \in L^2((0, T); H^3(I; \hmr{\mC}) \cap H^1_0(I; \hmr{\mC}))$. Define $\hmr{\ff} = \bT(f) \in L^1((0, T); \cD(\bA^*)')$. 
Then $\Phi \in C([0, T]; \bH)$ is a unique weak solution to the system 
\be 
\left\{\begin{array}{cl}
i \Phi_t = - \Delta \Phi  - \lambda \Phi   + \hmr{\ff}  & \mbox{ in } (0, T)  \times I, \\[6pt] 
\Phi(t, 0) = \Phi(t, 1) = 0 &  \mbox{ in } (0, T), \\[6pt]
\Phi(0, \cdot) = \Phi_0 & \mbox{ in } I,    
\end{array} \right. 
\ee 
(\hmr{in the sense given in \Cref{lem-WP}}) if and only if $\Phi \in C([0, T]; \bH)$ is a (weak) solution of the system 
\be 
\left\{\begin{array}{cl}
i \Phi_t = - \Delta \Phi  - \lambda \Phi   + f  & \mbox{ in } (0, T)  \times I, \\[6pt] 
\Phi(t, 0) = \Phi(t, 1) = 0 &  \mbox{ in } (0, T), \\[6pt]
\Phi(0, \cdot) = \Phi_0 & \mbox{ in } I,  
\end{array} \right. 
\ee 
in the sense that  
\be \label{def-WS-*}
i \frac{d}{dt} \langle \Phi, \varphi_k \rangle_{L^2(I)} = -  \langle \Delta \Phi, \varphi_k \rangle_{L^2(I)} -  \lambda \langle \Phi, \varphi_k \rangle_{L^2(I)}  + \langle f, \varphi_k \rangle_{L^2(I)} \mbox{ in } (0, T)
\ee
in the distributional sense for all $k \ge 1$. 
\end{lemma}

\begin{proof} Let $\varphi \in H^3(I; \mC) \cap H^1_0(I; \mC)$. We have 
$$
(1 + \lambda_k + \lambda_k^2 + \lambda_k^3) \langle \varphi, \varphi_k \rangle_{L^2(I)} = \langle \varphi, \varphi_k \rangle_{H^3(I)}  -  \varphi_{xx}(1) \hmr{\overline{ \varphi_{k,xxx} (1)}} + \varphi_{xx}(0) \hmr{\overline{\varphi_{\hmr{k}, xxx} (0)}}. 
$$
One can thus rewrite \hmr{\eqref{lem-WP-ak0}} under the form 
\be \label{lem-WP-L2}
i \frac{d}{dt} \langle \Phi, \varphi_k \rangle_{L^2(I)} = -  \langle \Phi, \Delta \varphi_k \rangle_{L^2(I)} + \langle f, \varphi_k \rangle_{L^2(I)} \mbox{ in } (0, T). 
\ee
The conclusion follows in the case $\lambda=0$. The general case follows similarly.  
\end{proof}

We next make a connection with the definition of weak solutions used in \cite{BL10}.  Let $e^{i \Delta t}: L^2(I; \mC) \to L^2(I; \mC)$ be defined by, for $\varphi \in L^2(I; \mC)$,   
\be
e^{i t \Delta} \varphi = \sum_{k=1}^\infty  \langle \varphi, \varphi_k \rangle_{L^2(I)} e^{- i \lambda_k t} \varphi_k,  
\ee
and,  for $\gamma \in \mR$, let  $e^{i (\Delta + \gamma) t}: L^2(I; \mC) \to L^2(I; \mC)$ be defined by, for $\varphi \in L^2(I; \mC)$, 
\be
e^{i t (\Delta + \gamma)} \varphi = \sum_{k=1}^\infty  \langle \varphi, \varphi_k \rangle_{L^2(I)} e^{ i ( - \lambda_k + \gamma) t} \varphi_k.   
\ee

We have
\begin{lemma} Let $0< T < T_0$ and  $\lambda \in \mR$. Let $\Phi_0 \in \bH$ and let $f \in L^2((0, T); H^3(I; \mC) \cap H^1_0(I; \mC))$. Define $\hmr{\ff} = \bT(f) \in L^1((0, T); \cD(\bA^*)')$. Then $\Phi \in C([0, T]; \bH)$ is a weak solution of \eqref{lem-WP-sys} if and only if 
\be \label{rem-WS}
\Phi(t, \cdot) = e^{i t (\Delta + \gamma)} \Phi_0 - i \int_0^t e^{i (t - s) (\Delta + \gamma)} \big(f(s, \cdot) + (\gamma-\lambda) \Phi(s, \cdot) \big) \, ds \mbox{ for } t \in [0, T].   
\ee
\end{lemma}

\begin{proof}
Assume that \eqref{rem-WS} holds for $\Phi \in C([0, T]; \bH)$. By taking the scalar product in $L^2(I)$ of the corresponding identity with $\varphi_k$, one derives that $\Phi$ is a weak solution of \eqref{lem-WP-sys} by \Cref{lem-product}. 

We next assume that $\Phi \in C([0, T]; \bH)$ is a weak solution of \eqref{lem-WP-sys}. We will prove \eqref{rem-WS}. Since 
$\Phi \in C([0, T]; \bH)$ is a weak solution of \eqref{lem-WP-sys}, we deduce from \Cref{lem-product} that
\begin{multline} \label{rem-WS-1-1}
\langle \Phi(t, \cdot), \varphi_k \rangle_{L^2(I)} = \langle e^{i t (\Delta + \gamma)} \Phi_0, \varphi_k \rangle_{L^2(I)} \\[6pt]
- \langle i \int_0^t e^{i (t - s) (\Delta + \lambda)}  \big(f(s, \cdot) + (\gamma-\lambda) \Phi(s, \cdot) \big) \, ds, \varphi_k \rangle_{L^2(I)} \mbox{ for } t \in [0, T]. 
\end{multline}
Set 
$$
\Psi(t, \cdot)= \int_0^t e^{i (t - s) (\Delta + \lambda)}  \big(f(s, \cdot) + (\gamma-\lambda) \Phi(s, \cdot) \big)  \, ds. 
$$
Since the space spanned by set of $(\varphi_k)_{k \ge 1}$ is dense in $L^2(I)$
and $\Psi \in C([0, T]; \bH)$, we obtain \eqref{rem-WS}. 
\end{proof}

\begin{remark} \rm In \cite[Proposition 2]{BL10}, the definition of the weak solutions in the sense of \eqref{rem-WS} is considered with $\gamma = 0$. 
\end{remark}

We next establish the well-posedness and stability of linear feedback systems.  It is convenient to denote
\be \label{def-cX}
\cX = \cX_T: = C([0, T]; \bH)
\ee
and to equip this space with the following standard norm 
\be \label{cX-norm}
\| \Phi \|_{\cX} : = \sup_{[0, T]} \|\Phi(t, \cdot) \|_{\bH}. 
\ee
Then $\cX$ is a Banach space. 

\hmb{\begin{lemma}\label{lem-WPWP-LN} Let $0< T \le T_0$, $\alpha > 0$, $\Psi \in C([0, T]; H^3(I; \mC) \cap H^1_0 (I; \mC))$, and
let $\bR_{ij} \in \cL (\bH_{rel})$ for $i, j = 1, 2$
be such that $\|\bR_{ij} \varphi_k\|_{C^3(\bar I)} \le \alpha k^3$ for some $\alpha > 0$ and for all $k \ge 1$, and $\| \Psi\|_{\cX} \le \alpha$. Denote 
$$
\bR \varphi = (\bR_{11} \varphi_1 + \bR_{12} \varphi_2) + i (\bR_{21} \varphi_1 + \bR_{22} \varphi_2) \mbox{ for } \varphi \in \bH, 
$$
where $\varphi_1 = \Re \varphi $ and $\varphi_2 = \Im \varphi$. Let $\Phi_0 \in \bH$ and $f \in L^2((0, T); H^3(I; \mC) \cap H^1_0(I; \mC))$.  There exists a unique weak solution $\Phi \in C([0, T]; \bH)$ to the system \footnote{The weak solution is understood in the sense given in \Cref{lem-WP}.}
\be
\left\{\begin{array}{cl}
i \Phi_t = - \Delta \Phi  - \lambda_1 \Phi + \bg & \mbox{ in } (0, T) \times I, \\[6pt] 
\Phi(t, 0) = \Phi(t, 1) = 0 &  \mbox{ in } (0, T), \\[6pt]
\Phi(0, \cdot) = \Phi_0 & \mbox{ in } I,   
\end{array} \right. 
\ee 
where $\bg \in L^1((0, T); \cD(A^*)')$ is defined by $\bg = \bT(g)$
with 
$$
g (t, \cdot) = u \Psi(t, \cdot) + f(t, \cdot), 
$$
where $u \in L^2((0, T); \mR)$ is such that $u= - B^*(\bR \Phi)$.   
Moreover, there exist a positive constant $C$ depending only on $T_0$ and $\alpha$  such that
\be 
\|  \Phi \|_{\cX} + \| u \|_{L^2(0, T)} \le C \Big( \| \Phi_0\|_{H^3(I)}  +  \| f\|_{L^2((0, T); H^3(I)} \Big) 
\mbox{ in } [0, T]. 
\ee
\end{lemma}}

\begin{proof} \hmb{Consider the map $\cK: L^2((0, T); \mR) \to L^2((0, T); \mR)$ defined by 
$$
\cK(u) = - B^*(\bR \Phi), 
$$
where $\Phi \in \cX$ is the unique solution of the system  
\be
\left\{\begin{array}{cl}
i \Phi_t = - \Delta \Phi  - \lambda_1 \Phi + \bg & \mbox{ in } (0, T) \times I, \\[6pt] 
\Phi(t, 0) = \Phi(t, 1) = 0 &  \mbox{ in } (0, T), \\[6pt]
\Phi(0, \cdot) = \Phi_0 & \mbox{ in } I.    
\end{array} \right. 
\ee 
Here $\bg \in L^1((0, T); \cD(A^*)')$ is defined by $\bg = \bT(g)$
with 
$$
g (t, \cdot) = u \Psi(t, \cdot) + f(t, \cdot).  
$$}

\hmb{For $u_j \in L^2((0, T); \mR) $ with $j=1, 2$, denote $\Phi_j$ the corresponding solutions. 
We have, by \Cref{lem-WP},  
$$
\|\cK(u_2) - \cK(u_1)\|_{L^2(0, T)} \le C_T \| u_2 - u_1 \|_{L^2(0, T)}. 
$$
Since one can choose $C_T$ small for $T$ small, the conclusion follows for small $T$. The conclusion in the general case follows from this case by dividing the interval $[0, T]$ into small sub-intervals.}
\end{proof}

We next study the local well-posedness and the stability of nonlinear feedback systems.

\hmb{\begin{lemma} \label{lem-WPWP-NL} Let $0< T \le T_0$, $\alpha > 0$, $\Psi \in C([0, T]; H^3(I; \mC) \cap H^1_0 (I; \mC))$, and 
and
let $\bR_{ij} \in \cL (\bH_{rel})$ for $i, j = 1, 2$
be such that $\|\bR_{ij} \varphi_k\|_{C^3(\bar I)} \le \alpha k^3$ for some $\alpha > 0$ and for all $k \ge 1$, and $\| \Psi\|_{\cX} \le \alpha$. Denote 
$$
\bR \varphi = (\bR_{11} \varphi_1 + \bR_{12} \varphi_2) + i (\bR_{21} \varphi_1 + \bR_{22} \varphi_2) \mbox{ for } \varphi \in \bH, 
$$
where $\varphi_1 = \Re \varphi $ and $\varphi_2 = \Im \varphi$. There exists a positive constant  $C$ depending only on $T_0$ and $\alpha$,  such that for $\Phi_0 \in \bH$ with 
\be \label{lem-WPWP-NL-small}
\| \Phi_0\|_{H^3(I)} \le \eps, 
\ee
there exists a unique solution $\Phi \in \cX$ to the system 
\be\label{lem-WPWP-NL-sys}
\left\{\begin{array}{cl}
i \Phi_t = - \Delta \Phi  - \lambda_1 \Phi + 
\bg & \mbox{ in } (0, T) \times I, \\[6pt] 
\Phi(t, 0) = \Phi(t, 1) = 0 &  \mbox{ in } (0, T), \\[6pt]
\Phi(0, \cdot) = \Phi_0 & \mbox{ in } I,  
\end{array} \right. 
\ee 
where $\bg \in L^1((0, T); \cD(A^*)')$ is defined by $\bg = \bT(g)$
with 
$$
g (t, x) = u (t) \Psi(t, x) + u (t) \mu (x) \Phi(t, x) \mbox{ for } (t, x) \in (0, T) \times (0, 1),  
$$
where $u \in L^2((0, T); \mR)$ is such that $u= - B^*(\bR \Phi)$.   
Moreover, there exist a positive constant $C$ depending only on $T_0$ and $\alpha$  such that
\be 
\|  \Phi(t, \cdot)\|_{H^3(I)} + \| u \|_{L^2((0, T); \mR)} \le C  \| \Phi_0\|_{H^3(I)}  \mbox{ in } [0, T]. 
\ee
\end{lemma}}

\begin{proof} \hmb{Let $\bPhi \in \cX$ be the solution of the system 
\be
\left\{\begin{array}{cl}
i \bPhi_t = - \Delta \bPhi  - \lambda_1 \bPhi + 
\bh & \mbox{ in } (0, T) \times I, \\[6pt] 
\bPhi(t, 0) = \bPhi(t, 1) = 0 &  \mbox{ in } (0, T), \\[6pt]
\bPhi(0, \cdot) = \bPhi_0 & \mbox{ in } I,  
\end{array} \right. 
\ee 
where $\bh \in L^1((0, T); \cD(A^*)')$ is defined by $\bh = \bT(h)$
with 
$$
h (t, \cdot) = u \Psi(t, \cdot), 
$$
where $u \in L^2((0, T); \mR)$ is such that $u= - B^*(\bR \bPhi)$. By \Cref{lem-WPWP-LN}, one has 
$$
\| \bPhi \|_{\cX} \le C \eps \mbox{ provided that \eqref{lem-WPWP-NL-small} holds}.  
$$}

\hmb{One now considers the map $\cK: \overline{B_{\cX}(\bPhi, \eps)} \to \cX$ defined by  $\cK(\wPhi)$, for $\wPhi \in  \overline{B_{\cX}(\bPhi, \eps)}$ (the closed ball centered at $\bPhi$ and of radius $\eps$ in $\cX$),  is the solution $\Phi \in \cX$ of the system 
\be\label{lem-WP-CS-NL-sys}
\left\{\begin{array}{cl}
i \Phi_t = - \Delta \Phi  - \lambda_1 \Phi + 
\bg & \mbox{ in } (0, T) \times I, \\[6pt] 
\Phi(t, 0) = \Phi(t, 1) = 0 &  \mbox{ in } (0, T), \\[6pt]
\Phi(0, \cdot) = \Phi_0 & \mbox{ in } I,  
\end{array} \right. 
\ee 
where $\bg \in L^1((0, T); \cD(A^*)')$ is defined by $\bg = \bT(g)$
with 
$$
g (t, \cdot) = u \Psi(t, \cdot) + u \mu \wPhi(t, \cdot), 
$$
where $u \in L^2((0, T); \mR)$ is such that $u= - B^*(\bR \Phi)$. Then, by \Cref{lem-WPWP-LN},  
$$
\| \Phi \|_{\cX} + \| B^*(\bR \Phi) \|_{L^2(0, T)} \le C \eps. 
$$
Using this fact, we derive from \Cref{lem-WPWP-LN} that 
$$
\| \Phi - \bPhi \|_{\cX} \le  C \| B^*(\bR \Phi) \mu \wPhi \|_{L^2((0, T); H^3(I))} \le C \eps^2. 
$$
Thus, for $\eps$ sufficiently small, one has
$$
\cK: \overline{B_{\cX}(\bPhi, \eps)}  \to  \overline{B_{\cX}(\bPhi, \eps)}. 
$$}

\hmb{We next show that $\cK$ is contracting.  Let $\wPhi_j \in \cap \overline{B(\bPhi, \eps)}$ for $j=1, 2$. Denote  $\Phi_j$ the corresponding solution and $u_j$ the corresponding control 
for $j=1, 2$. We have, by \Cref{lem-WPWP-LN},  
\be \label{lem-WPWP-NL-p1}
\|\Phi_2 - \Phi_1\|_{\cX} \le C\| u_2 \wPhi_2 - u_1 \wPhi_1 \|_{L^2((0, T); \mH^3(I))}. 
\ee
We have 
\begin{multline*}
\| u_2 \wPhi_2 - u_1 \wPhi_1 \|_{L^2((0, T); \mH^3(I))} \le C \Big( \| (u_2 - u_1) \wPhi_2  \|_{L^2((0, T); \mH^3(I))}  + \| u_1( \wPhi_2 -  \wPhi_1) \|_{L^2((0, T); \mH^3(I))} \Big) \\[6pt]
\le C \eps \Big( \|u_2 - u_1 \|_{L^2(0, T)} + \|\wPhi_2 - \wPhi_1 \|_{\cX} \Big) \mathop{\le}^{\Cref{lem-WPWP-LN}} C \eps \Big(\| u_2 \wPhi_2 - u_1 \wPhi_1 \|_{L^2((0, T); \mH^3)} + \|\wPhi_2 - \wPhi_1 \|_{\cX} \Big). 
\end{multline*}
It follows that, for $\eps$ sufficiently small,  
\be\label{lem-WPWP-NL-p2}
\| u_2 \wPhi_2 - u_1 \wPhi_1 \|_{L^2((0, T); \mH^3(I))}  \le C \eps \|\wPhi_2 - \wPhi_1 \|_{\cX}. 
\ee}

\hmb{Combining \eqref{lem-WPWP-NL-p1} and \eqref{lem-WPWP-NL-p2} yields 
$$
\|\Phi_2 - \Phi_1\|_{\cX} \le C \eps \| \wPhi_2 - \wPhi_1 \|_{\cX}.   
$$}

\hmb{Therefore, for $\eps$ sufficiently small, $\cK$ is contracting and the conclusion follows.} 
\end{proof}

We next translate the previous well-posedness result to the semi-group related to $A$ defined in \eqref{def-A}, which involves the definition of $Q$ and the feedback. 
We only do it for \Cref{lem-WP}. The statement and the proof of the corresponding variants of  \hmr{\Cref{lem-WPWP-LN} and \Cref{lem-WPWP-NL}} are omitted to avoid repetition. To this end, we first introduce the following definition. 

\begin{definition} \label{def-TT} Let $T>0$ and $\lambda \in \mR$, and let $(\cA, \cD(\cA))$ be defined by \eqref{def-cA}.  Define 
$$
\TT: L^2( (0, T); H^3(I; \mR^2) \cap H^1_0 (I; \mR^2)) \to L^1((0, T); \cD(\cA^*)')
$$
by 
\be\label{lem-WP-g1}
\langle \TT(g)(t, \cdot), \varphi \rangle_{\cD(\cA^*)', \cD(\cA^*)} = \langle g(t, \cdot), \varphi \rangle_{H^3(I)}  -  \langle g(t, 1),  \varphi_{xxx} (1) \rangle_{\mR^2} + \langle g(t, 0),  \varphi_{xxx} (0) \rangle_{\mR^2}.  
\ee
\end{definition}

Concerning \Cref{lem-WP}, we have the following result. 

\begin{lemma} \label{lem-WP-A}
Let $0 < T < T_0$, $\lambda \in \mR$, and let $\cA: \cD(\cA) \subset \mH \to \mH$ be defined by 
\be \label{def-cA}
\cA y = 
\left(\begin{array}{cc} - \Delta y_2 - \lambda y_2 \\[6pt]
\Delta y_1 + \lambda y_1
\end{array}\right)
 \quad \mbox{ and } \quad \cD(\cA) = \Big\{ y \in \mH; \cA y \in \mH \Big\}.  
\ee
Let $y_0 \in \mH$ and $\bg \in L^1((0, T); \cD(\cA^*)')$.  There exists a unique weak solution $y \in C([0, T]; \cD(\cA^*)')$ to the system 
\be\label{lem-WP-sys-A}
\left\{\begin{array}{cl}
y\hmr{'} = \cA y + \bg & \mbox{ in } (0, T)  \times I, \\[6pt] 
y(0, \cdot) = y_0 & \mbox{ in } I. 
\end{array} \right. 
\ee 
Let $g = (g_1, g_2)\tr \in L^2((0, T); H^3(I; \hmr{\mR^2}) \cap H^1_0(I; \hmr{\mR^2}))$. Define $\bg  = \TT(g) \in L^1((0, T); \cD(\cA^*)')$. Then $y \in C([0, T]; \mH)$. Moreover, $y = (y_1, y_2)\tr$ is a weak solution of \eqref{lem-WP-sys-A} if and only if $\Phi: = y_1 + i y_2$ is a weak solution of \eqref{lem-WP-sys} with 
$$
\Phi_0  = y_1(0, \cdot) + i y_2(0, \cdot) \quad \mbox{ and } f = -  g_2 + i g_1,  
$$
and $\ff = \bT(f)$. We also have
$$
\|y(t)\|_{\mH} \le C\Big( \| y_0\|_{\mH}  + \| g\|_{L^2((0, T); H^3(I))} \Big) 
\mbox{ in } [0, T], 
$$
where $C$ is a positive constant depending only on $T_0$.  
We also have
$$
\|y(t)\|_{\mH} \le C\Big( \| y_0\|_{\mH}  + \| g\|_{L^2((0, T); H^3(I))} \Big) 
\mbox{ in } [0, T], 
$$
\hmb{Let $\cR = (\cR_1, \cR_2) \tr \in \cL(\mH)$ with 
$\|\cR(\Phi^{1}_k)\|_{C^3(\bar I)}, \|\cR(\Phi^{2}_k)\|_{C^3(\bar I)} \le \alpha k^3$ for some $\alpha > 0$ and for all $k \ge 1$. Then
\be
B^* (\cR y) \in L^2((0, T); \mR^2) \quad \mbox{ and } \quad  \| B^* y \|_{L^2(0, T)} \le C_1 \| y_0\|_{\mH}  + C_2 \| g\|_{L^2((0, T); H^3(I))},  
\ee
for some positive constants $C_1, C_2$ depending only on $T_0$ and $\alpha$. Moreover, $C_2$ can be chosen such that $C_{2} \to 0$ as $T_0 \to 0$. }
\end{lemma}

\begin{proof}[Proof of \Cref{lem-WP-A}] By \Cref{pro-WP} in the appendix, there exists a unique weak solution $y \in C([0, T]; \cD(\cA^*)')$ \hmr{of \eqref{lem-WP-sys-A}}. Let $\Phi \in C([0, T]; \bH)$ be the unique weak solution of \eqref{lem-WP-sys} with 
$$
\Phi_0  = y_1(0, \cdot) + i y_2(0, \cdot) \quad \mbox{ and } \quad f = -  g_2 + i g_1. 
$$
Let $y_1$ and $y_2$ be the real part and the imaginary part of $\Phi$, respectively,  and denote $y = (y_1, y_2) \tr$. Then 
$$
y \in C([0, T]; \mH). 
$$
By \Cref{lem-WP}, it suffices to prove that $y$ is a weak solution of \eqref{lem-WP-sys-A}. This follows from the definition of weak solutions associated with $\bA$. 
\end{proof}

\section{Rapid stabilization - Proof of \Cref{thm1-S} and \Cref{thm2-S}} \label{sect-Rapid}

This section containing two subsections is devoted to the proof of \Cref{thm1-S} and \Cref{thm2-S}. The proof of \Cref{thm1-S} is given in the first subsection and the proof of \Cref{thm2-S} is given in the second one. 

\subsection{Proof of \Cref{thm1-S}}

Denote 
$$
A_\lambda = A + \lambda I.
$$ 
Let $(y_\lambda, \ty_{\lambda})\tr \in C([0, +\infty); \mH)$ be the unique weak solution of the system 
\be \label{thm-main1-sys-hytylambda}
\left\{\begin{array}{c}
y_\lambda' = A_\lambda y_\lambda - B B^* \ty_\lambda \quad \mbox{ in } (0, T), \\[6pt]
\ty_\lambda'  = - A_\lambda^* \ty_\lambda  \quad \mbox{ in } (0, T), \\[6pt]
\hy_\lambda(0) = y_0, \quad  \ty_\lambda(0) = \cQ^{-1} y_0. 
\end{array} 
\right. 
\ee

Let  $\tau \in (0, T]$ and $\varphi_\tau \in \mH$, and let $\varphi \in C([0, \tau]; \mH)$ be the unique weak solution of 
\be \label{thm-main1-varphi} 
\left\{\begin{array}{c}
\varphi ' = - A_\lambda^* \varphi  \mbox{ in } (0, \tau), \\[6pt]
\varphi(\tau) = \varphi_\tau.  
\end{array} \right. 
\ee
Applying \cite[Lemma \hmr{3.1}]{Ng-Riccati} for $A_\lambda$ with $t = \tau$, we derive from \eqref{thm-main1-sys-hytylambda} and \eqref{thm-main1-varphi} that  
\be \label{thm-main1-p2}
 \langle y_\lambda (\tau), \varphi(\tau) \rangle_{\mH} - \langle y_\lambda (0), \varphi(0) \rangle_{\mH} =  - \int_0^\tau \langle B^* \ty_\lambda (s), B^* \varphi(s) \rangle_{\mR} \, ds. 
\ee
Applying \cite[Lemma \hmr{4.1}]{Ng-Riccati} to $\ty_\lambda(\tau - \cdot)$ and $\varphi (\tau - \cdot)$, we obtain \hmr{from \eqref{thm-main1-sys-hytylambda} and \eqref{thm-main1-varphi} that}
\be\label{thm-main1-p3}
\langle Q \ty_\lambda(0), \varphi(0) \rangle_{\mH} - \langle Q \ty_\lambda(\tau), \varphi (\tau) \rangle_{\mH}  \\[6pt]
 =  \int_0^\tau \langle B^*\ty_\lambda(s), B^*\varphi (s) \rangle_{\mR}  \, ds.    
\ee
Summing \eqref{thm-main1-p2} and \eqref{thm-main1-p3}, after using the fact that $Q \ty_\lambda(0) = y_\lambda (0)$,  we deduce that  
$$
\langle y_\lambda(\tau) - Q \ty_\lambda(\tau), \varphi(\tau) \rangle_{\mH}  = 0. 
$$

Since $\varphi(\tau) \in \mH$ is arbitrary, we derive that  
\be\label{thm-main1-p1-000}
y_\lambda(\tau) - Q \ty_\lambda(\tau) =  0. 
\ee

Set 
$$
\hy(t) = e^{-\lambda t} y_{\lambda}(t), \quad \mbox{ and } \quad \ty (t)= e^{-\lambda t} \ty_{\lambda}(t). 
$$
Then, from \eqref{thm-main1-sys-hytylambda} and \eqref{thm-main1-p1-000},  we have 
\be \label{thm-main1-sys-hyty}
\left\{\begin{array}{c}
\hy' = A \hy - B B^* \ty \quad \mbox{ in } (0, T), \\[6pt]
\ty'  = - A^* \ty  -2 \lambda \ty \quad \mbox{ in } (0, T), \\[6pt]
\hy(0) = y_0, \quad  \ty(0) = \cQ^{-1} y_0, 
\end{array} 
\right. 
\ee
and 
\be \label{thm-main1-p1-000-z}
\hy - Q \ty =  0 \mbox{ in } [0, T].  
\ee

\hmr{From \eqref{thm-main1-sys-hyty}, using the fact that $A$ is skew-symmetric and the definition of $B$ and $B^*$,  we have 
\be
\frac{d}{dt} \langle \hy, \Phi_1 \rangle_{\mH} = 0 \mbox{ in } (0, T). 
\ee
Recall that $\Phi_1$ is defined in \eqref{def-F}. Since $y_0 \in \mH_{1, \sharp}$, we derive that  $\hy(t) \in \mH_{1, \sharp}$ in $[0, T]$.}  It follows that 
$$
y = \hy \mbox{ in } [0, T]. 
$$

Since
$$
\ty' = A \ty - 2 \lambda \ty, 
$$
and $A$ is skew-adjoint by \Cref{lem-A},  it follows that 
\be\label{thm-main1-p1-222}
\| \ty(t)\|_{\mH} = e^{-2 \lambda t} \| \ty_0 \|_{\mH}.  
\ee

The conclusion now follows from \eqref{thm-main1-p1-222} and the fact that $y(t) = Q 
\ty(t)$ for $t \ge 0$. \qed

\subsection{Proof of \Cref{thm2-S}}
Let $y_P(t)$ be the projection of $y(t)$ into $\mH_{1, \sharp}$ using the $\mH$-scalar product and set $z(t) = y(t) - y_P(t)$. Let $b :  [0, + \infty) \to \mR$ be the real function such that 
\be
z =  b \Phi_1
\ee
(recall that $\Phi_1 = (\varphi_1, 0)\tr$). Note that 
\be \label{thm2-S-def-u}
u(t) = - B^* \cQ^{-1} y_P (t). 
\ee
Since $A z = 0$, we derive from \eqref{thm2-S-sys} that 
\be \label{thm2-S-p1}
y' = y_P' + z'  = A y_P +  B u    +  u F(y - \Phi_1) \mbox{ for } t > 0.  
\ee
Taking the scalar product in $[L^2(I)]^2$ of  this equation with $\Phi_1$ and integrating by parts, we obtain, by \Cref{lem-product},  
\be \label{thm2-S-def-b'}
b' = u \langle F(y - \Phi_1), \Phi_1 \rangle_{L^2(I)}. 
\ee
We derive from \eqref{thm2-S-p1} and \eqref{thm2-S-def-b'} that 
\be \label{thm2-S-eq-Y}
y_P'  = A y_P +  B u + u F(y - \Phi_1) - z' =   A y_P +  B u + u F(y - \Phi_1) - b' \Phi_1.  
\ee

Fix $T_0 > 0$ and let $0 < T < T_0$. \hmr{We claim that 
\be\label{estimate-tY-1}
\|y_P(t)\|_{\mH} \le \| \cQ^{-1} \|_{\cL(\mH_{1, \sharp})} e^{-2 \lambda t} \| y_P(0)\|_{\mH} + C_{\lambda, T_0} \|u F(y - \Phi_1) - b' \Phi_1\|_{L^2 ((0, T);  H^3(I))} \mbox{ in } [0, T].  
\ee
Indeed, we have 
\be \label{estimate-tY-1-1}
y_{P} = y_{P, 1} + y_{P, 2} \quad \mbox{ in } [0, T], 
\ee
where, with $u_{1} = - B^* \cQ^{-1} y_{P, 1}$, 
$$
\left\{\begin{array}{cl} 
y_{P, 1}'  = A y_{P, 1} +  B u_{1} \mbox{ in } (0, T), \\[6pt]
y_{P, 1}(0) = y_P(0), 
\end{array} \right.
$$
and
$$
\left\{\begin{array}{cl}y_{P, 2}'  = A y_{P, 2} + B (u - u_1) + u F(y - \Phi_1) - b' \Phi_1 \mbox{ in } (0, T), \\[6pt]
y_{P, 2}(0) = 0. 
\end{array} \right.
$$
Applying \Cref{thm1-S} to $y_{P, 1}$ and  \Cref{lem-WPWP-LN} (its variant in the spirit of \Cref{lem-WP-A}) to $y_{P, 2}
$ after noting $y_{P, 2} = y_{P} - y_{P_1}$ and $u - u_1 = -  B^* \cQ^{-1} (y_P - y_{P_1}) = -  B^* \cQ^{-1} y_{P, 2}$, we obtain 
\be\label{estimate-tY-1-2}
\|y_{P, 1}(t)\|_{\mH} \le \| \cQ^{-1} \|_{\cL(\mH_{1, \sharp})} e^{-2 \lambda t} \| y_P(0)\|_{\mH}  \mbox{ in } [0, T], 
\ee
and  
\be\label{estimate-tY-1-3}
\|y_{P, 2}(t)\|_{\mH} \le  C_{\lambda, T_0} \|u F(y - \Phi_1) - b' \Phi_1\|_{L^2 ((0, T);  H^3(I))}. 
\ee
Claim \eqref{estimate-tY-1} now follows from \eqref{estimate-tY-1-1}, \eqref{estimate-tY-1-2}, and \eqref{estimate-tY-1-3}. }

Since, for $t \in [0, T]$, 
$$
\| u\|_{L^2 (I)} +  \|y - \Phi_1 \|_{L^\infty ((0, T);  \mH)} \mathop{\le}^{\hmr{\Cref{lem-WPWP-NL}}} C_{T_0} \hmr{\| y(0) - \Phi_1 \|_{\mH}},
$$
$$
\|F(y - \Phi_1) \|_{L^\infty ((0, T);  H^3(I))} \mathop{\le}^{\eqref{def-F}} C_{T_0} \|y - \Phi_1 \|_{L^\infty ((0, T);  \mH)},
$$
and  
$$
\| b'\|_{L^2 (0, T)} \mathop{\le}^{\eqref{thm2-S-def-b'}} C_{T_0} \| u\|_{L^2(0, T)} \| y - \Phi_1 \|_{L^\infty(0, T)} \le C_{T_0} \hmr{\|y(0) - \Phi_1 \|_{\mH}^2}, 
$$
it follows from \eqref{estimate-tY-1} that 
\be \label{estimate-tY}
\|y_P(t)\|_{\mH} \le  \| \cQ^{-1} \|_{\cL(\mH_{1, \sharp})} e^{-2 \lambda t} \| \hmr{y_0 - \Phi_1}\|_{\mH} + C_{\lambda, T_0}\hmr{ \|y(0) - \Phi_1\|_{\mH}^2}\mbox{ in } [0, T].
\ee

By the conservation of $L^2$-norm of $y$,  we derive from the fact that $z$ and $y_P$ are orthogonal in $L^2(I; \mR^2)$, 
$$
b(t)^2 + \|y_P(t) \|^2_{L^2(I)} = 1. 
$$
It follows that  if $\|y(t) - \Phi_1 \|_{L^2(I)} \le 1/4$ in $[0, T]$, then 
\be\label{thm1-S-pp1}
\| b(t) - 1\|_{L^2(I)} \le C \|y_P(t) \|^2_{L^2(I)}. 
\ee

By taking $\eps_0$ sufficiently small and $\| y_0 \hmr{- \Phi_1} \|_{\mH} \le \eps_0$, we derive from \eqref{estimate-tY}, and \eqref{thm1-S-pp1} that 
\be
\|y(t) - \Phi_1\|_{\mH} \le    2 \| \cQ^{-1} \|_{\cL(\mH_{1, \sharp})} e^{-2 \lambda t} \| y(0) - \Phi_1\|_{\mH} \mbox{ in } [0, T_0].   
\ee
Taking $T_0$ large enough such that $2 \| \cQ^{-1} \|_{\cL(\mH_{1, \sharp})} e^{-2 \lambda T_0} \le e^{- 2 \hlambda T_0}$, we then can repeat the argument for the interval $[T_0, 2 T_0]$, \dots, and obtain the conclusion. 

\section{Finite-time stabilization - Proof of \Cref{thm-FT-LN}} \label{sect-Finite}

This section containing two subsections is devoted to the proof of \Cref{thm-FT-LN}. In the first subsection, we establish the cost of controls for the linearized system in small time. The proof of \Cref{thm-FT-LN} using the results in the first section is given in the second one. 

\subsection{Cost of control of the linearized system for small time}

We begin this section by establishing an upper bound of the cost of control of the linearized system for small time. 

\begin{proposition} \label{pro-cost-upper}  Let $\mu \in H^3(I, \mR)$ be such that \eqref{cond-mu} holds. For all $\Psi_0 \in \bH_{1, \sharp}$, there exists $u \in L^2(0, T); \mR)$ such that 
$$
\Psi(T, \cdot) = \Psi_0
$$
and
$$
\| u\|_{L^2(0, T)} \le e^{\frac{C}{T}} \| \Psi_0\|_{\bH},  
$$
where $\Psi \in C([0, T]; \bH)$ is the unique weak solution of the system 
\be
\left\{\begin{array}{cl}
i \Psi_t = - \Delta \Psi  - \lambda_1 \Psi - \bT(u(t) \mu  \varphi_1)  & \mbox{ in } (0, T) \times I, \\[6pt] 
\Psi(t, 0) = \Psi(t, 1) = 0 &  \mbox{ in } (0, T), \\[6pt]
\Psi(0, \cdot) = 0 & \mbox{ in } I. 
\end{array} \right. 
\ee 
\end{proposition}

Recall that $\bT$ is defined in \Cref{def-T}.

\begin{proof}  By a translation of time, it suffices to prove the following result. For all $\Psi_0 \in \bH_{1, \sharp}$, there exists $u \in L^2((-T/2, T/2); \mR)$ such that 
$$
\Psi(T/2, \cdot) = \Psi_0
$$
and
$$
\| u\|_{L^2(-T/2, T/2)} \le e^{\frac{C}{T}} \| \Psi_0\|_{\bH},  
$$
where $\Psi \in C([-T/2, T/2]; \bH)$ is the unique weak solution of the system 
\be
\left\{\begin{array}{cl}
i \Psi_t = - \Delta \Psi  - \lambda_1 \Psi - \bT(u(t) \mu  \varphi_1)  & \mbox{ in } (-T/2, T/2) \times I, \\[6pt] 
\Psi(t, 0) = \Psi(t, 1) = 0 &  \mbox{ in } (-T/2, T/2), \\[6pt]
\Psi(-T/2, \cdot) = 0 & \mbox{ in } I. 
\end{array} \right. 
\ee 

The proof of this fact is based on the moment method, see, e.g.,  \cite{TT07}.  We represent $\Psi$ under the form
$$
\Psi(t, x) = \sum_{k \ge 1} a_k (t) \varphi_k (x) \mbox{ in } (-T/2, T/2) \times I. 
$$
We then have, see the proof of \Cref{lem-product},  
$$
i a_k' = (\lambda_k - \lambda_1) a_k - c_k u(t) \mbox{ in } (-T/2, T/2),   
$$
where 
\be \label{pro-cost-upper-ck}
c_k =  \langle \mu \varphi_1, \varphi_k \rangle_{L^2(I)} \mbox{ in } (-T/2, T/2). 
\ee
Thus
$$
a_k' = - i (\lambda_k - \lambda_1) a_k + i c_k u(t) \mbox{ in } (-T/2, T/2).   
$$
Since $a_k(-T/2) = 0$ for $k \ge 1$, we then have 
$$
a_k (T/2) =   i c_k \int_{-T/2}^{T/2} e^{- i (\lambda_k - \lambda_1) (T/2-s)} u (s) \, ds.    
$$
Set
\be \label{pro-cost-wk}
v (t) = u(T/2 - t), \quad \omega_k = \lambda_k - \lambda_1 \mbox{ for } k \ge 1, 
\ee
and
\be \label{pro-cost-dk0}
d_k = a_k(T/2)/ (i c_k) \mbox{ for } k \ge 1. 
\ee
Since $a_1(T/2) = \Im a_1(T/2)$ for $\Phi_0 \in \bH_{1, \sharp}$, it follows from \eqref{pro-cost-upper-ck} that 
$$
d_1 \in \mR. 
$$ 
We then have
$$
\int_{-T/2}^{T/2} e^{- i \omega_k s} v (s) \, ds  = d_k \mbox{ in } (-T/2, T/2) \mbox{ for } k \ge 1,      
$$
 
By \eqref{cond-mu} and \eqref{pro-cost-dk0}, we have 
\be \label{pro-cost-dk}
\sum_{k \ge 1} |d_k|^2 \le C \|\Psi_0 \|_{\bH}^2.  
\ee

By \cite[Lemma 4.1]{TT07} there exists $\beta > 0$ such that, for all $k \ge 1$, \footnote{\cite[Lemma 4.1]{TT07} only gives the estimate for the first term; nevertheless, the estimate for the second term can be done in the same manner.} 
\be \label{pro-cost-t1}
\prod_{n \ge 1, n \neq k} \left|1 - \frac{z}{\omega_n - \omega_k} \right| + \prod_{n \ge 1, n \neq k} \left|1 + \frac{z}{\omega_n + \omega_k} \right|  \le \beta e^{\beta |z|^{1/2}} \mbox{ for } z \in \mC.  
\ee
By \cite[Lemma 4.2]{TT07}, for all $\gamma > 1$, there exist $C = C(\gamma)$ (independent of $T \in (0, T_0)$) and  an analytic function $H$ such that 
\be \label{pro-cost-t2}
H(0) = 1 \quad \mbox{ and } \quad |H(z)| \le e^{\frac{C}{T}} e^{- \gamma |z|^{1/2}} \mbox{ for } z \in \mR, \quad |H(z)| \le C_T e^{T |\Im z |/4}. 
\ee

Fix $\gamma > 2 \beta$ and a corresponding analytic function $H$. For $N 
\ge 2$, we define the function $\xi_N: \mC \to \mC$ as follows, for $z \in \mC$,  
\begin{multline}
\xi_N(z) = \sum_{k = 2}^N d_k H(z - \omega_k) \prod_{n \ge 1, n \neq k} \Big(1 - \frac{z - \omega_k}{\omega_n - \omega_k} \Big) \prod_{l \ge 1}   \Big(1 - \frac{z - \omega_k}{-\omega_l - \omega_k} \Big)  \\[6pt]
+ \sum_{k = 2}^N \bar{d_k} H(- z - \omega_k) \prod_{n \ge 1, n \neq k} \Big(1 - \frac{z + \omega_k}{-\omega_n + \omega_k} \Big) \prod_{l \ge 1}   \Big(1 - \frac{z + \omega_k}{\omega_l + \omega_k} \Big) \\[6pt]
+ d_1 H(z) \prod_{n \ge 2} \Big(1 - \frac{z}{\omega_n} \Big) \prod_{l \ge 2}   \Big(1 - \frac{z}{-\omega_l} \Big).
\end{multline}
It follows from \eqref{pro-cost-dk}, \eqref{pro-cost-t1}, and \eqref{pro-cost-t2} that the function $\xi_N$ is well-defined and is analytic on $\mC$. From the definition of $\xi$, we have
\be \label{pro-cost-xi1}
\xi_N(\omega_k) = d_k, \quad \xi_N(-\omega_k) = \bar{d_k} \mbox{ for } 2 \le k \le N, 
\ee
and  
\be \label{pro-cost-xi2}
\xi_N(0) = \xi_N(\omega_1)  = d_1. 
\ee

For all $c>0$, there exists $c_1 > 0$ such that it holds  
\begin{multline*}
\int_{\mR} e^{-c|z - \omega_m|^{1/2}} e^{-c|z - \omega_n|^{1/2}} \, dz + \int_{\mR} e^{-c|z + \omega_m|^{1/2}} e^{-c|z + \omega_n|^{1/2}} \, dz  \\[6pt] 
+ \int_{\mR} e^{-c|z - \omega_m|^{1/2}} e^{-c|z + \omega_n|^{1/2}} \, dz \le c_1 e^{- \frac{c}{2} |\omega_{\hmr{m}} - \omega_n|^{1/2}}. 
\end{multline*}
We derive from \eqref{pro-cost-dk},  \eqref{pro-cost-t1}, and \eqref{pro-cost-t2} that 
\be \label{pro-cost-L2}
\|\xi_N\|_{L^2(\mR)} \le e^{\frac{C}{T}} \|\Psi_0\|_{\bH}. 
\ee
and 
$$
\mbox{the restriction of } \xi_N \mbox{ on $\mR$ is a Cauchy sequence in $L^2(\mR)$}. 
$$
By Paley-Wiener's theorem, see, e.g., \cite[Theorem 19.3]{Rudin-RC}, there thus exists $v_N \in L^2(-T/2, T/2; \mC)$ such that $\hv_N = \xi_N$ and $\hmr{v_N}$ is a Cauchy sequence in $L^2(-T/2, T/2)$. Let $v$ be the limit of the sequence $(v_N)$ in $L^2(-T/2, T/2)$.  Set 
$$
u = \frac{1}{2} \Big(v + \overline{v} \Big).
$$
Then, for $k \ge 1$,  
$$
\hu (\omega_k) = \frac{1}{2} \lim_{N \to + \infty}\Big(\xi_N (\omega_k) + \overline{\xi_N} (-\omega_k)  \Big) \mathop{=}^{\eqref{pro-cost-xi1}, \eqref{pro-cost-xi2}} d_k. 
$$
Here we used the fact that $d_1$ is real.  The conclusion follows since 
$$
\| u\|_{L^2(\mR)} \le C \| v\|_{L^2(\mR)} \le C \limsup_{N \to + \infty} \| \xi_N\|_{L^2(\mR)} \mathop{\le}^{\eqref{pro-cost-L2}} e^{\frac{C}{T}} \|\Psi_0\|_{\bH}
$$
The proof is complete. 
\end{proof}

For the completeness, we next establish a lower bound of the cost of control of the linearized system for small time. To this end, we first prove the following result. 

\begin{proposition}\label{pro-S-M} Let $0 < \eps < 1/2$, $\eps^3 < T < 3/2$.  Let $\mu \in L^1((0, T); 
\mR)$ be such that $\int_{I} \mu \varphi_1^2 \, dx \neq 0$.  If  $u \in L^1((0, T); 
\mC)$ is a control which steers the control system 
\begin{equation} \left\{
\begin{array}{c}
i v_t = - \eps v_{xx} - \eps \lambda_1 v + \frac{i}{4 \eps} v - \eps u(t) \mu(x) \varphi_1   \mbox{ in } (0, T) \times (0, 1), \\[6pt]
v(t, 0) = v(t, 1) \mbox{ in } (0, T)
\end{array} \right. 
\end{equation} 
from  $\varphi_1 $ at time $0$ to $0$ at the time $T$ in the sense that there exists $v \in L^2((0, T); H^1_0(I)) 
\cap C([0, T]; L^2(I))$ such that 
\begin{multline}
i \frac{d}{dt} \langle v, \varphi_k \rangle_{L^2(I)} = -  \eps \langle v, \Delta \varphi_k \rangle_{L^2(I)} -  
\eps \lambda_1 \langle v, \varphi_k \rangle_{L^2(I)}  \\[6pt]
+ \frac{i}{4 \eps} \langle v, \varphi_k \rangle_{L^2(I)}
- \eps u(t) \langle \mu \varphi_1, \varphi_k \rangle_{L^2(I)} \mbox{ in } (0, T)
\end{multline}
in the distributional sense for all $k 
\ge 1$, and $v(0, \cdot) = \varphi_1$ in $I$ and $v(T, 
\cdot) = 0$ in $I$, then
\begin{equation*}
\ln  \| u \|_{L^1(0, T)}   \ge \frac{1}{\eps} \Big(\frac{1}{2} - \frac{T}{4}\Big) -  C \ln \eps^{-1},  
\end{equation*}
for some positive constant $C$ independent of $\eps$ and $T$.
\end{proposition}

\begin{proof} The proof uses tools from complex analysis, see, e.g., \cite{Koosis09}. Define
\begin{equation*}
\Lambda_k :=   \eps \lambda_k - \eps \lambda_1 + \frac{i}{4 \eps}
 \quad \mbox{ and } \quad \Phi_k(t, x) = \varphi_k (x) e^{- i \Lambda_k t} \mbox{ in } \mR_+ \times (0,1).
\end{equation*}
One can check that 
\begin{equation}
\left\{\begin{array}{cl}
\dsp i \Phi_{k, t} = -  \eps \Phi_{k, xx} - \eps \lambda_1 \Phi_k + \frac{i}{4 \eps} \Phi_k = 0& \mbox{ in } \mR_+ \times [0, 1], \\[6pt]
\Phi_k(t, 0) = \Phi_k(t, 1) = 0  & \mbox{ for } t \in  \mR_+.
\end{array}\right.
\end{equation}
Multiplying the equation of $v$ by $\Phi_k$, integrating by parts, and using the fact $v(T, \cdot) =0$,  we have
\begin{equation}\label{dual1}
i \int_0^1 v(0, x) \Phi_k(0, x) \, dx  = \eps \int_0^T u(t) e^{- i \Lambda_k t}  \, dt
\int_I \mu \varphi_1 \varphi_k \, d x. 
\end{equation}
Define, for $z \in \mC$,
\begin{equation}\label{def-FFF}
\cF(z) :=   \int_{-T/2}^{T/2} u(t + T/2) e^{- i z t} \, dt = e^{i T z /2} \int_0^T u(t) e^{ - i z t} \, dt.
\end{equation}
It follows from \eqref{dual1} that, with $c_1 = \int_I \mu \varphi_1^2 \, dx$, 
\begin{equation}\label{zero1}
\cF (  \Lambda_1 ) = \frac{i}{c_1 \eps } e^{i \Lambda_1 T /2}  \quad 
\mbox{ and } \quad \cF ( \Lambda_k ) =0  \mbox{ for } k \ge 2.
\end{equation}
Applying  the representation of entire functions of exponential type for $\cF$, see e.g.,  \cite[page 56]{Koosis09}, we derive from \eqref{zero1} that, for $z  \in \mC$ with $\Im z > 0$,
\begin{equation}\label{complex-analysis} 
\ln |\cF (z)| \le I_0(z)  + I_1(z) + \sigma \Im( z),
\end{equation}
where
\begin{equation}
I_0 (z)= \sum_{k \ge 2} \ln \frac{\big|\Lambda_k - z \big|}{\big|
\bar \Lambda_k - z \big|}, \quad I_1(z) = \frac{\Im(z)}{\pi} \int_{-\infty}^\infty \frac{\ln |\cF(\tau)|}{|\tau - z|^2} \, d\tau,
\end{equation}
and
\begin{equation}
\sigma = \limsup_{y \to + \infty}  \frac{\ln |\cF( i y)|}{y}.
\end{equation}
From the definition of $\cF$ in \eqref{def-FFF}, we have, for $y \in \mR_+$,  
\begin{equation*}
|\cF(i y)| \le \| u\|_{L^1(0, T)} e^{T y/2}. 
\end{equation*}
This implies 
\begin{equation}\label{est-F-1}
\sigma = \limsup_{y \to + \infty}  \frac{\ln |\cF(i y)|}{y} \le T/2. 
\end{equation}
From  \eqref{zero1},  we derive that 
\begin{equation}\label{Fb1}
\ln |F(\Lambda_1)| \ge - \frac{T}{8 \eps}  +   C \ln \eps^{-1}.
\end{equation}
Here and in what follows in this proof, $C$ denotes a positive constant independent of $k$ and it can change from one place to another. 
Similar to  \cite[(2.597)]{Coron07} (see also \cite{CG05}),  we obtain 
\begin{multline}\label{est-I0-1}
I_0 (\Lambda_{1} ) =  \sum_{k \ge 2} \ln \frac{\eps (k^2-1) \pi^2}{ \big([\eps (k^2 - 1) \pi^2]^2 + [1/(2\eps)]^2 \big)^{1/2}} \le  \sum_{k \ge 2} \ln \frac{ \eps^2 k^2 \pi^2}{ \big([\eps^2 k^2  \pi^2]^2 + [1/2]^{2} \big)^{1/2}} \\[6pt] 
\le  \int_1^\infty \ln \left(\frac{\eps^2 \pi^2  x^2}{ (\eps^4 \pi^4  x^4 + 1/4 \big)^{1/2}} \right) \, dx = \frac{1}{\eps \pi \sqrt{2}} \int_{\eps \pi  \sqrt{2} }^\infty \ln \Big(\frac{x^2}{\sqrt{x^4 + 1}} \Big) \, dx. 
\end{multline}
 Since 
$$
\int_{0}^\infty \ln \Big(\frac{x^2}{\sqrt{x^4 + 1}} \Big) \, dx = - \frac{ \pi} { \sqrt{2}}, 
$$
it follows from \eqref{est-I0-1} that 
\begin{equation}\label{I0}
I_0 (\Lambda_{1} )  \le - \frac{1}{2 \eps}  + C \ln \eps^{-1}.
\end{equation}
We next estimate $I_1 (\Lambda_{1} )$.  From \eqref{def-FFF},  we have, for $s \in \mR$, 
$$
\ln |\cF(s)|  \le \| u\|_{L^1(0, T)},   
$$
which yields 
\begin{equation}\label{I1-1}
I_1(\Lambda_{1}) \le \frac{1}{4 \eps \pi} \int_{-\infty}^\infty \frac{\ln \|u \|_{L^1(0, T)} }{s^2 + (1/4\eps)^2} \, ds. 
\end{equation}
Since, for $a>0$,
\begin{equation*}
\frac{a}{\pi} \int_{-\infty}^\infty \frac{1}{s^2 +a^2} \, ds = 1,
\end{equation*}
it follows from \eqref{I1-1} that
\begin{equation}\label{I1-2}
I_1 \big(\Lambda_{1} \big) \le \ln  \|u \|_{L^1(0, T)}.
\end{equation}
Combining \eqref{complex-analysis}, \eqref{Fb1}, \eqref{I0}, and \eqref{I1-2}  yields
\begin{equation*}
-\frac{T}{8\eps} \le - \frac{1}{2\eps } + \ln \| u \|_{L^1(0, T)}  + \frac{T}{8 \eps} + C \ln \eps^{-1}. 
\end{equation*}
This implies
\begin{equation*}
\ln  \| u \|_{L^1(0, T)}   \ge \frac{1}{\eps} \Big(\frac{1}{2} - \frac{T}{4}\Big) -  C \ln \eps^{-1}. 
\end{equation*}
The proof is complete.  
\end{proof}

We are ready to obtain a lower bound for the cost of control viewing \Cref{lem-product}.

\begin{proposition} \label{pro-cost-lower} Let $\mu \in L^1((0, T); 
\mC)$ be such that $\int_{I} \mu \varphi_1^2 \, dx \neq 0$, and let $0 < T < 1$. Then if $u \in L^1((0, T); \mC)$ is a control which steers the control system 
\begin{equation} \left\{
\begin{array}{c}
i v_t = - v_{xx} -  \lambda_1 v  -  u(t) \mu(x) \varphi_1   \mbox{ in } (0, T) \times (0, 1), \\[6pt]
v(t, 0) = v(t, 1) \mbox{ in } (0, T)
\end{array} \right. 
\end{equation} 
from  $\varphi_1$ at time $0$ to $0$ at the time $T$ in the sense that there exists $v \in L^2((0, T); H^1_0(I)) 
\cap C([0, T]; L^2(I))$ such that 
\be
i \frac{d}{dt} \langle v, \varphi_k \rangle_{L^2(I)} = -   \langle v, \Delta \varphi_k \rangle_{L^2(I)} -  
 \lambda_1 \langle v, \varphi_k \rangle_{L^2(I)}  
- u(t) \langle \mu \varphi_1, \varphi_k \rangle_{L^2(I)} \mbox{ in } (0, T)
\ee
in the distributional sense for all $k 
\ge 1$, and $v(0, \cdot) = 
\varphi_1$ in $I$ and $v(T, 
\cdot) = 0$ in $I$, then 
\begin{equation*}
\ln \| u\|_{L^1(0, T)}  \ge  \frac{1}{4 T} - C  \ln T^{-1},
\end{equation*}
for some positive constant $C$ independent of $T$. 
\end{proposition}

\begin{proof} Define
\begin{equation*}
 \tv(t, x) : = v(\eps t, x)  e^{-\frac{t}{ 4 \eps}} \mbox{ for } (t, x) \in (0, T / \eps) \times (0, 1)
\end{equation*}
and set $\tu(t) = u(\eps t) e^{-\frac{t}{ 4 \eps}}$ for $t \in (0, T / \eps)$.
Then
\begin{equation*}
i \tv_t  = - \eps \tv_{xx} - \eps \lambda_1 \tv + \frac{i}{4 \eps} \tv -  \eps \tu (t) \mu (x) \varphi_1 \mbox{ in } (0, T/\eps) \times (0, 1). 
\end{equation*}
Applying \Cref{pro-S-M} to $\tv$ with $(\eps, T)  = (T,T)$, we have
\begin{equation*}
\ln \| u\|_{L^1(0, T)}  \ge   \ln \| \tu \|_{L^1(0, 1)}  \ge \frac{1}{T} \Big(\frac{1}{2} - \frac{1}{4}\Big) -  C \ln T^{-1} \ge \frac{1}{4 T} - C  \ln T^{-1},
\end{equation*}
which is the conclusion. 
\end{proof}

\begin{remark} \rm
Similar arguments as in the proof of \Cref{pro-cost-lower} can be found in \cite{Lissy12}. For the boundary controls, the cost of controls for small time is also of the order $e^{C/T}$, see, e.g., \cite{TT07}. 
\end{remark}

\begin{remark} \rm It is shown in \cite{Ng22-Cost} that the cost of controls of the heat equation depends on the support of the data and the controlled region. This is based on the strategy of Lebeau and Robbiano \cite{LR95} and the three-sphere inequalities with partial data established by Nguyen \cite{Ng-CALR-O-M}. It would be interesting to study whether or not the cost of controls depends on the support of the initial data for the KdV system. 
\end{remark}

\subsection{Proof of \Cref{thm-FT-LN}}

We first give an estimate for $Q = Q(\lambda)$. The following result is a direct consequence of \Cref{pro-cost-upper} and Hilbert uniqueness method. 

\begin{proposition} \label{pro-cost} Let $\mu \in H^3(I; \mR)$ be such that \eqref{cond-mu} holds and  let $0< T < T_0$.  We have, for some positive constant $C$ independent of $T$, 
$$
\int_0^T | B^* e^{-s A^*} z|^2  \ge e^{- \frac{C}{T}} \|z\|_{\mH}^2 \mbox{ for all } z \in \mH_{1, \sharp}.  
$$
\hmr{Consequently, for  $\lambda \ge  \lambda_0$ and for $Q = Q(\lambda)$ defined by \eqref{def-Q}, it holds 
$$
\| Q \| \le C, 
$$
for some positive constant $C$ independent of $\lambda$.}
\end{proposition}

Using \Cref{pro-cost}, we can prove the following result. 

\begin{lemma} \label{lem-Qlambda} Let $\lambda \ge  \lambda_0$ and let $Q = Q(\lambda)$ be defined by \eqref{def-Q}. There exists a positive constant $C$ independent of $\lambda$ such that 
\be
\langle Q z, z \rangle_{\mH} \ge e^{-C \sqrt{\lambda}} \| z\|_{\mH}^2 \mbox{ for all } z \in \mH_{1, \sharp}. 
\ee \label{lem-Qlambda-cl1}
\hmr{Consequently, 
\be \label{lem-Qlambda-cl2}
\| \cQ^{-1} \| \le e^{C \sqrt{\lambda}},  
\ee
where $\cQ$ is defined in \eqref{def-cQ}.} 
\end{lemma}

\begin{proof}
We have 
\begin{multline}
\langle Q z, z \rangle_{\mH} = \int_{0}^\infty e^{-2 \lambda s} | B^* e^{-s A^*} z|^2 \, ds \ge \int_{1/\sqrt{\lambda}}^{2/  \sqrt{\lambda}} e^{-2 \lambda s} | B^* e^{-s A^*} z|^2 \, ds \\[6pt]
\ge \frac{1}{\sqrt{\lambda}}  e^{-4 \sqrt{\lambda} } \int_{1/\sqrt{\lambda}}^{2/\sqrt{\lambda}}| B^* e^{-s A^*} z|^2 \, ds  
\mathop{\ge}^{\Cref{pro-cost}}  \frac{C_1}{\sqrt{\lambda}}  e^{-4 \sqrt{\lambda} } e^{- C_2\sqrt{\lambda}} \| e^{- \frac{1}{\sqrt{\lambda}} A^*} z\|_{\mH}, 
\end{multline}
which yields 
\be
\langle Q z, z \rangle_{\mH} \ge  \frac{C_1}{\sqrt{\lambda}}  e^{-4 \sqrt{\lambda} } e^{- C_2 \sqrt{\lambda}} \| z\|_{\mH}. 
\ee
\hmr{Assertion \eqref{lem-Qlambda-cl1} follows. Assertion \eqref{lem-Qlambda-cl2} is just a consequence of \eqref{lem-Qlambda-cl1} and the definition of $\cQ$.} 
\end{proof}

\Cref{thm-FT-LN} is now a consequence of the following result.

\begin{proposition} \label{pro-FN-LN}  Let $\mu \in H^3(I, \mR)$ be such that \eqref{cond-mu} holds and  let $T > 0$.  Let $(t_n)$ be an increasing sequence that converges to $T$ with $t_0 = 0$ and let $(\lambda_n) \subset \mR_+$ be an increasing sequence. Define, for $t_{n} \le t < t_{n+1}$ and $n \ge 0$, 
$$
\cK(t, z) = - B^* Q_{n}^{-1} \proj_{\mH_{1, \sharp}} z \mbox{ for } z \in \mH, 
$$
where $Q_n = Q(\lambda_n)$ defined by \eqref{def-Q} with $\lambda = \lambda_n$. 
Set $s_0 = 0$ and $s_n = \sum_{k=0}^{n-1} \lambda_k(t_{k+1} - t_k)$ for $n \ge 1$. Let
$y \in C([0, T); \mH)$ \hmr{be the unique solution} of system \eqref{sys-LN1-AB} with 
$$
u(t) = \cK(t, y(t, \cdot)) \mbox{ for } t \in [0, T). 
$$
There exists a positive constant $\gamma$ such that, if for large $n$, 
\be  \label{pro-FN-LN-as1}
(t_{n+1} - t_n) \lambda_n \ge \gamma \sqrt{\lambda_n}, 
\ee
then it holds, for $t_{n-1} \le t \le t_{n}$ and for $n \ge 1$,  
\be \label{pro-FN-LN-cl1}
\|y(t, \cdot) \|_{\mH} \le e^{ - s_{n-1} + C n} \| y_0\|_{\mH}
\ee
and 
\be \label{pro-FN-LN-cl2}
|u(t)| \le C e^{-s_{n-1}/4 + C n} \| y_0 \|_{\mH}, 
\ee
for some positive constant $C$ independent of $n$. In particular, if, in addition, we have that
\be  \label{pro-FN-LN-as2}
\lim_{n \to + \infty} \frac{s_n}{n + \sqrt{\lambda_{n+1}}} = + \infty, 
\ee
then 
\be \label{pro-FN-LN-cl3}
y(t, \cdot) \to 0 \mbox{ in } \mH \mbox{ as } t \to T_{-}
\ee
and 
\be \label{pro-FN-LN-cl4}
u(t, \cdot) \to 0 \mbox{ as } t \to T_{-}. 
\ee
\end{proposition}

\begin{remark} \rm There are sequences $(t_n)$ and $(\lambda_n)$ which satisfy the conditions given in the above proposition, for example,  $t_n = T - T/n^2$ and $\lambda_n = n^8$ for large $n$. 
\end{remark}

\begin{proof} \hmr{Applying \Cref{thm1-S}, we have, for $n \ge 1$,  
\be \label{pro-FN-LN-p1}
\| \cQ^{-1} y(t) \|_{\mH}  = e^{-2 \lambda_n(t - t_{n-1})} \| \cQ^{-1} y(t_{n-1}) \|_{\mH} \mbox{ for } t \in [t_{n-1}, t_n]. 
\ee
We derive from \Cref{lem-Qlambda} that, for $n \ge 1$,  
\be
\| y(t) \|_{\mH} \le e^{- 2 \lambda_{n-1} (t_n - t_{n-1}) + C \sqrt{\lambda_{n-1}}} \| y(t_{n-1}) \|_{\mH} \mbox{ for } t \in [t_{n-1}, t_n]. 
\ee
In particular, one has 
\be
\| y(t_n) \|_{\mH} \le e^{- 2 \lambda_{n-1} (t_n - t_{n-1}) + C \sqrt{\lambda_{n-1}}} \| y(t_{n-1}) \|_{\mH} \mbox{ for } n \ge 1.
\ee
It follows that 
\be\label{pro-FN-LN-p2}
\| y(t_{n}) \|_{\mH} \le e^{- 2 s_{n-1} + C \sum_{k=0}^{n-1} \sqrt{\lambda_k}} \| y_0 \|_{\mH} \mbox{ for } n \ge 1. 
\ee
Using \eqref{pro-FN-LN-as1}, we derive \eqref{pro-FN-LN-as1} from \eqref{pro-FN-LN-p1} and \eqref{pro-FN-LN-p2}.}

We have, by \eqref{def-B*}, for $t_{n-1} \le t \le t_n$ and for $n \ge 1$,  
$$
|u(t)| = |B^*Q_{n-1}^{-1} y(t, \cdot)| \le C e^{C \sqrt{\lambda_{n-1}}} \| y(t, \cdot) \|_{\mH}. 
$$
\hmr{Assertion \eqref{pro-FN-LN-cl2} now follows from from \eqref{pro-FN-LN-p1} and \eqref{pro-FN-LN-p2} after using \eqref{pro-FN-LN-as1}.}
\end{proof}

\begin{remark}  \rm We are not able to extend the finite-time stabilization to the nonlinear setting. This is due to the fact \hmr{that} we cannot ensure the well-posedness for the time interval $[t_n, t_{n+1})$ for large $n$. 
\end{remark}

\appendix

\section{Control systems associated with operator semi-groups} \label{sect-appendixA}

In this section, we recall and establish some facts on the control systems associated with a strongly continuous semigroup. The standard references are \cite{Zwart96,EN00,Coron06,BDDM07,TW09}.

Let $\cH$ and $\cU$ be two Hilbert spaces which denote the state space and the control space, respectively. The corresponding scalar products are  $\langle \cdot, \cdot \rangle_{\cH}$ and $\langle \cdot, \cdot \rangle_{\cU}$, and the corresponding norms are $\| \cdot \|_{\cH}$ and $\| \cdot \|_{\cU}$. Let $\big(S(t) \big)_{t \ge 0} \subset \cL(\cH)$ be  a strongly continuous  semi-group on $\cH$. Let $(\cA, \cD(\cA))$ be the infinitesimal generator of $\big(S(t) \big)_{t \ge 0}$ and denote $S(t)^*$ the adjoint of $S(t)$ for $t \ge 0$. Then $\big( S(t)^* \big)_{t \ge 0}$ is also a strongly continuous semigroup of continuous linear operators and its infinitesimal generator is 
$(\cA^*, \cD(\cA^*))$, which is the adjoint of $(\cA, \cD(\cA))$. As usual, we equip the domain $\cD(\cA^*)$ with the scalar product 
$$
\langle z_1, z_2 \rangle_{\cD(\cA^*)} = \langle z_1, z_2 \rangle_{\cH} +  \langle \cA^* z_1, \cA^* z_2 \rangle_{\cH} \mbox{ for } z_1, z_2 \in \cD(\cA^*).  
$$
Then $\cD(\cA^*)$ is a Hilbert space. Denote $\cD(\cA^*)'$ the dual space of $\cD(\cA^*)$ with respect to $\cH$. Then 
$$
\cD(\cA^*) \subset \cH \subset \cD(\cA^*)'. 
$$ 
Let 
$$
\cB \in \cL(\cU, \cD(A^*)'). 
$$

As usual, we equip the domain $\cD(A^*)$ with the scalar product 
$$
\langle z_1, z_2 \rangle_{\cD(\cA^*)} = \langle z_1, z_2 \rangle_{\cH} +  \langle \cA^* z_1, \cA^* z_2 \rangle_{\cH} \mbox{ for } z_1, z_2 \in \cD(\cA^*).  
$$
Then $\cD(A^*)$ is a Hilbert space. Denote $\cD(\cA^*)'$ the dual space of $\cD(\cA^*)$ with respect to $\cH$. Then 
$$
\cD(\cA^*) \subset \cH \subset \cD(\cA^*)'. 
$$ 
Let 
$$
\cB \in \cL(\hmr{\cU}, \cD(\cA^*)'). 
$$
Consider the control system  
\be \label{CS-G}
\left\{\begin{array}{c}
y' = \cA y + f + \cB u  \mbox{ in } t \in (0, T), \\[6pt]
y(0) = y_0, 
\end{array} \right. 
\ee
with $y_0 \in \cD(\cA^*)'$, and $f \in L^1((0, T); (\cD(\cA^*))')$ and $u \in L^1((0,T); \cU)$. We are interested in weak solutions of \eqref{CS-G}.  

\begin{definition} \label{def-WS-cA}
A weak solution $y$ of \eqref{CS-G} is understood as an element $y \in C([0, T]; \big(\cD(\cA^*) \big)')$ such that 
\be\label{meaning-CS-G}
\left\{\begin{array}{c} \frac{d}{dt} \langle y, \varphi \rangle_{\cH}  = \langle \cA y + f + \cB u, \varphi \rangle_{\cH} \mbox{ in } (0, T) \\[6pt]
y(0) =  y_0
\end{array}\right. \mbox{ for all } \varphi \in \cD({\cA^*}^\infty)
\ee
for which 
\begin{itemize}
\item[$i)$] the differential equation in \eqref{meaning-CS-G} is understood in the distributional sense, 
\item[$ii)$] the term $ \langle \cA y + f + \cB u, \varphi \rangle_{\cH} $ is understood as $\langle y, \cA^*\varphi \rangle_{\cH}  +  \langle f, \varphi \rangle_{\cD(\cA^*)', \cD(\cA^*)} + \langle u, \cB^*\varphi \rangle_{\cU}$. 
\end{itemize}
\end{definition}

The convention in $ii)$ will be used from later on. Recall that $\cD({\cA^*}^\infty)$ is dense in $\cD(\cA^*)$, see e.g.,  \cite[Proposition 1.7]{EN00}. The following result is on the well-posedness of weak solutions of \eqref{CS-G}. 

\begin{proposition}\label{pro-WP} Let $T>0$, $y_0 \in \cD (\cA^*)'$,  $u \in L^1((0, T); \cU)$,  and $f \in L^1((0, T); \cD (\cA^*)')$.  Then   $y \in C([0, T], \cD(\cA^*)')$ is a weak solution of \eqref{CS-G}
if and only if, with $\tf := f + \cB u$, it holds \footnote{This identity is understood in $\cD(A^*)'$, i.e., $\langle y(t), \varphi \rangle_{\mH}  = \langle S(t) y_0, \varphi \rangle_{\mH} + \int_0^t \langle S(t - s) \tf(s), \varphi \rangle_{\mH} \, ds$ in  $[0, T]$ for all $\varphi \in \cD({\cA^*}^\infty)$. The solutions defined by \eqref{pro-WP-identity} are called mild solutions.}
\be\label{pro-WP-identity}
y(t) = S(t) y_0 + \int_0^t S(t - s) \tf(s) \, ds \mbox{ for } t \in [0, T]. 
\ee
\end{proposition}

\begin{proof} We first prove that $y \in C([0, T]; \cD(\cA^*)')$ is a weak solution of \eqref{CS-G} if and only if $y \in C([0, T]; \cD(\cA^*)')$ and \eqref{pro-WP-identity} holds. 

Assume first that $y \in C([0, T]; \cD(\cA^*)')$ and \eqref{pro-WP-identity} holds. 
We will prove that $y$ is a weak solution of \eqref{CS-G}. Here and in what follows, for notational ease, we denote $\langle \cdot, \cdot \rangle_{\cH}$
by $\langle \cdot, \cdot \rangle$. From \eqref{pro-WP-identity}, we obtain, with $\varphi \in \cD({A^*}^\infty)$,   
$$
\langle y(t), \varphi \rangle =  \langle S(t) y_0 + \int_0^t S(t - s)\tf(s) \, ds - y_0, \varphi \rangle \mbox{ for } t \in [0, T]. 
$$

Set 
\be
\psi(t) =  \langle S(t) y_0 + \int_0^t S(t - s) \tf(s) \, ds - y_0, \varphi \rangle \mbox{ for } t \in [0, T]. 
\ee
Then, for $t \in [0, T]$,  
\be \label{psi-pro-WP}
\psi(t) =  \langle y_0, S(t)^*\varphi \rangle +  \int_0^t \langle \tf(s), S(t - s)^* \varphi \rangle \, ds - \langle y_0, \varphi \rangle. 
\ee
Since $\varphi \in \cD({\cA^*}^\infty)$, we derive from \eqref{psi-pro-WP} that, for $t \in [0, T]$,   
\begin{multline*}
\psi'(t) =  \langle y_0, \cA^*S(t)^* \varphi \rangle +  \int_0^t \langle \tf(s), \cA^*S(t - s)^* \varphi \rangle \, ds \rangle + \langle \tf(t), \varphi \rangle \\[6pt]
= \langle y_0, S(t)^* \cA^* \varphi \rangle +  \int_0^t \langle \tf(s), S(t - s)^* \cA^*\varphi \rangle \, ds + \langle \tf(t), \varphi \rangle,  
\end{multline*}
which yields, by \eqref{pro-WP-identity},  
\be \label{pro-WP-p1}
\psi'(t)  =  \langle y, \cA^* \varphi \rangle + \langle \tf(t), \varphi \rangle. 
\ee
Integrating \eqref{pro-WP-p1} and using \eqref{psi-pro-WP} and \eqref{pro-WP-identity}, we obtain, in $[0, T]$,  
$$
\langle y(t) - y_0, \varphi \rangle  = \int_0^t \langle y, \cA^* \varphi \rangle + \langle \tf(t), \varphi \rangle, 
$$
which in turn implies \eqref{meaning-CS-G}.

\medskip 
We now prove that if $y$ is a weak solution of \eqref{CS-G}, then $y$ satisfies \eqref{pro-WP-identity}. We first assume that $f \in C([0, T]; \cH)$ and $y_0 \in \cD(\cA)$.
For $t > 0$, set, with $\varphi \in \cD({A^*}^\infty)$ and $s \in [0, t]$,
$$
\chi(s) =  \langle S(t-s) y(s), \varphi \rangle  =  \langle y(s), S(t-s)^* \varphi \rangle. 
$$
Then
$$
\chi'(s) =  - \langle y(s), S(t-s)^*A^* \varphi  \rangle + \langle A y(s) + \tf(s), S(t-s)^* \varphi \rangle = \langle S(t-s) \tf(s), \varphi \rangle. 
$$
It follows that 
$$
\chi(t) - \chi(0) = \int_0^t \langle S(t-s) \tf(s), \varphi \rangle \, ds, 
$$
which yields the identity. The proof in the general case follows by density. 

The proof is complete. 
\end{proof}

\begin{remark} \rm  The equivalence between weak solutions and mild solutions was first proved in the case $B$ is bounded and $f \in C([0, T]; \cH)$ by Ball \cite{Ball77}, see also \cite[Chapter 1 of Part II]{BDDM07} for related results when $\cB$ is bounded. 
\end{remark}


\providecommand{\bysame}{\leavevmode\hbox to3em{\hrulefill}\thinspace}
\providecommand{\MR}{\relax\ifhmode\unskip\space\fi MR }
\providecommand{\MRhref}[2]{%
  \href{http://www.ams.org/mathscinet-getitem?mr=#1}{#2}
}
\providecommand{\href}[2]{#2}

\end{document}